\documentclass[12pt]{article}
\usepackage{amsmath,amssymb,color}
\usepackage{times}
\usepackage{array}

\usepackage{graphicx}
\input epsf

\def\qed{\hfill$\Box$}
\def\bR{\mathbb R}

\def\bT{\mathbb T}
\def\ep{\epsilon}

\def\sgn{\hbox{sgn}\,}

\def\cJ{{\mathcal J}}

\def\cL{{\mathcal L}}

\def\cP{{\cal P}}

\numberwithin{equation}{section}
\newtheorem{remark}{Remark}[section]
\newtheorem{proposition}{Proposition}[section]
\newtheorem{theorem}{Theorem}[section]
\newtheorem{definition}{Definition}[section]
\newtheorem{lemma}{Lemma}[section]
\newtheorem{corollary}{Corollary}[section]
\begin{document}


\title{\bf Oscillating facets}

\author{Milena Matusik\\ 
Institute of Mathematics\\ University of Gda\'{n}sk\\
ul. Wita Stwosza 57, 80-952 Gda\'{n}sk, POLAND\\
{\tt Milena.Matusik@mat.ug.edu.pl}\\
Piotr Rybka\\Faculty of Mathematics, Informatics and  Mechanics\\
The University of Warsaw\\
ul. Banacha 2, 02-097 Warsaw, POLAND\\{\tt rybka@mimuw.edu.pl}}

\maketitle
\date{}

\abstract{We study a singular one-dimensional parabolic problem with
initial
 data in the $BV$ space, the energy space, for various boundary data. We pay special attention to Dirichlet conditions, which need not satisfied in a pointwise manner. We study the facet creation process and the extinction of solutions
caused by the evolution of facets. Our major tool is the comparison principle
provided by the theory of viscosity solutions developed in \cite{miyory}.}

\bigskip
\noindent{\bf keywords:} strongly singular parabolic equations, facets, anisotropy, viscosity  solutions, extinction time


\section{Introduction}\label{sin}
We study here a sudden diffusion problem
\begin{equation} \label{rng} 
\frac{\partial u}{\partial t}= 
\left( {\cL}(u_x)\right) _x\quad\hbox{in }I_T:=(a,b)\times(0,T).
\end{equation}
In other words, we assume that a monotone function $\cL$ has at least one jump discontinuity. More precisely, we consider
\begin{equation}\label{eL}
 \cL(p) = \sgn (p+1) + \sgn (p-1),
\end{equation}
for various boundary conditions. We pay special attention to Dirichlet data, which need not be satisfied in a pointwise manner.

Problems like this one appear in the models of crystal growth. We  may regard
(\ref{rng}) as a Gibbs-Thomson law, involving weighted mean curvature for a two-dimensional crystal, written in
local coordinates. Another approach to crystal growth is presented by Spohn,
\cite{spohn}, who discusses an equation like (\ref{rng}), but his choice of $\cL$
involves also a degenerate term, which we drop for the sake of the simplicity of
analysis.

Problem (\ref{rng}) is interesting even if we consider $\cL(p) = \sgn (p)$,
see \cite{becen}, \cite{taylor},  \cite{fukui-giga},  \cite{giga-kob},
\cite{spohn}, \cite{mury12} and the references therein. The nonlinearity, we
consider here, appears naturally, when we consider a corner formed by two
evolving facets. By a facet we mean a part of the graph of a solution to
(\ref{rng}) with the slope corresponding to a jump in $\cL$. In the present
case, facets have  slope $\pm 1$. Facets will  be defined rigorously in
Subsection \ref{serg}.

In a series of papers, \cite{gr2006}, \cite{gr2008}, \cite{gr2009},  \cite{gigory}, \cite{gigory2}, we studied evolution of the so-called bent rectangles by the weighted mean curvature flow,
\begin{equation}\label{wmcf}
 \beta V = \kappa_\gamma + \sigma.
\end{equation}
The point is, the corners of these bent rectangles were formed by facets meeting at the right angle. If we choose the local coordinate system in a proper way, then after simplifications preserving the main difficulties, system (\ref{wmcf}) looks like (\ref{rng}), this is presented in \cite{gigory2}. The main point is that nonlinearity (\ref{eL}) supports facets with different slopes.

Our main objective is to study interactions of facets, especially in the case of
oscillating data. In order to make equation (\ref{rng}) well-posed, we augment it with initial condition 
\begin{equation} \label{rni}
u(x,0)= u_0(x)\quad\hbox{for } x\in I:=(a,b)
\end{equation}
and either Dirichlet,
\begin{equation} \label{rnD}
 u(a)=A, \qquad u(b)=B,
\end{equation}
periodic
\begin{equation} \label{rnP}
 u(a)=u(b)
\end{equation}
or Neumann boundary data,
\begin{equation} \label{rnN}
 \cL(u_x)|_{\partial I} = 0.
\end{equation}
We put a cap on the oscillatory behavior of the data by requiring that $u_0\in
BV$.

One of the emerging problems is the meaning of (\ref{rnD}). It turns out
that our definition of solutions is too weak to guarantee that (\ref{rnD}) is
satisfied in a pointwise manner: the trace of solutions to (\ref{rng}),
(\ref{rni}) and (\ref{rnD}) need not be equal to the boundary data. We
elaborate on this in Definition \ref{d1} in subsection \ref{sub2.1} and \S \ref{set}.

We notice that initial conditions from $BV$ may have infinitely many facets of
different slopes. We would like to determine if this is possible for any 
solution at $t>0$.
We shall see that most of the facet interactions are resolved instantly. Thus, at $t>0$, we may have only a finite number of facets with non-zero curvature, see Theorem \ref{tmf}.

Our task involves re-examining the existence result of \cite{mury}, because we
consider less regular data than there. It is helpful to
observe that  (\ref{rng}) is formally a gradient flow of functional $E$ on
$L^2(I)$ defined by 
$$
E(u) = \int_I W(u_x)\, ds,
$$
where $W(p) = |p+1 |+ |p-1|$. Obviously, $E$ is well-defined iff $u\in BV(I)$. Thus, we will seek solutions with finite energy if $u_0\in BV(I)$.

We also have to discuss the notion of a solution to (\ref{rng})
defined in \cite[Theorem 1]{mury}, because smooth solutions to the approximating
system satisfying the Dirichlet data need not satisfy them in the limit if the
convergence is too weak. In order to expose the issue of the Dirichlet
boundary data, we will present explicit solutions in Proposition \ref{pr-ex}. We
make additional comments when we characterize the steady states in Section
\ref{set}.

We mention in passing that by a solution we mean a pair $(u, \Omega)$, where
$\Omega(\cdot, t)$ is a selection of the subdifferential 
$\partial E(u(\cdot, t))$. 
More details will be given in
Section \ref{ser}. It turns out that studying values of $\Omega$ and its continuity gives a lot of information about solutions. In many instances, see Section \ref{set}, Subsection \ref{serg},
this is our major tool.

Once the existence of solutions is established, we will
characterize the steady states for all three boundary conditions. This is done
in Section \ref{set}. In principle, they belong to $BV$.
We will see  that if $A$ does not differ much from $B$, then the steady states
are Lipschitz continuous  functions satisfying (\ref{rnD}) and  such that 
\begin{equation} \label{rnr}
 |u_x| \le 1.
\end{equation}
In turns out however, that if the  difference $B-A$ is big, then there are also discontinuous steady states belonging to $BV$. In other words the $BV$ regularity of the steady states is optimal. 
On the other hand, we note that all  Lipschitz functions satisfying (\ref{rnr}) are the steady states with the homogeneous Neumann data. 

We see  the multitude of the steady states. Condition (\ref{rnr}) permits seemingly
unchecked oscillations. This seems surprising. We will present two
justifications of this phenomenon.
Namely, we notice at all times $t>0$ there are only finitely many facets with nonzero
curvature, see Theorem \ref{tmf}. The other explanation is that our solutions
are viscosity solutions in the sense of \cite{miyory}.  We will see that 
in Section \ref{lepkie}. In addition, the theory of viscosity
solutions gives us a powerful tool like the Comparison Principle, see Theorem
\ref{comprin}. It is used in the proof of the main result of subsection
\ref{sect}, i.e. estimates on the extinction time of solutions. 

In Section \ref{sep},
we study the regularizing action of the flow when $u_0\in BV$ and the extinction
times of solutions. 
We constructed solutions in the energy space, i.e., $u(t) \in BV$, by way of examples, see Section \ref{ssdisc}, we shall see that discontinuities in $u_0$ persist. A more
interesting observation is that $u_0\in BV (I)$ implies that $u_t\in L^2(I_T)$
and this
statement  carries a lot of information about regularity and oscillatory behavior of solutions.
Namely, for almost all $t>0$ we have $u_t(\cdot,t) \in L^2(I)$. This implies
that the number of facets with non-zero curvature is finite for almost all
$t>0$, see Theorem \ref{tmf}. The argument is based on the observation that $\Omega_x(\cdot, t)\in
\partial E(u(\cdot, t))$. 

As we mentioned, $u(t)$ may have jumps as well as $u_x(t)$. We will see, see Theorem \ref{proB}, that jumps of
the derivative may not be arbitrary. This fact is well-known for the crystalline
motion, see \cite{GG-kin} and the references therein.
If at $x_0$ the interval with endpoints $u^+_x(x_0)$, $u^-_x(x_0)$ contains any of the singular slopes from $\{-1,1\}$, then immediately the missing facet is created for $t>0$. A similar statement holds if $u$ has a jump discontinuity at $x_0$.

The conclusion that for almost all $t>0$ solution $u(\cdot,t)$ has a finite
number of non-zero curvature facets permits more detailed studies of the
equations of facet motion. In Section \ref{sep} we concentrate on estimates in
terms of
initial data. We see that facet interaction is the main mechanism for the
extinction of solutions.  Due to
the fact that we have only a finite number of
moving facets the task is easier. Our main tool is the comparison principle,
Theorem \ref{comprin}, for
viscosity solutions established in \cite{miyory}. In particular, we give in
Theorem \ref{T-ext} a simple (but not a
closed formula) estimate for the extinction time in terms of the data. This is
the content of Subsection \ref{sect}.


This paper is devoted to the study of (\ref{rng}) with homogeneous boundary
conditions. The reason is that the theory of equations like (\ref{rng}) with a
forcing term 
$$
u_t =\cL(u_x)_x +f
$$
needs further development. Presently, we have only partial results, see e.g.
\cite{chano}, \cite{meyer}, \cite{gigi}
\cite{gigory}, \cite{gigory2}, \cite{mucha}. This is why
we make the restriction on the data. 


\section{Existence reexamined}\label{ser}

We introduce the definition of solutions to (\ref{rng}) with either Dirichlet,
periodic or
Neumann boundary data. We stress that it is well-known, see
\cite{dirichlet.andreu}, that the case of Dirichlet is more difficult and the boundary
need not be satisfied pointwise. We note that the problem becomes more apparent when
we want to interpret (\ref{rng}) as a gradient flow. The obvious functional on $L^2(I)$,
$$
E(u) = \left\{
\begin{array}{ll}
 \int_{I} |u_x -1|+|u_x +1| & \hbox{for } u\in BV(I), \ \gamma u(a)=A, \ \gamma u(b)=B,\\
+\infty & \hbox{else}
\end{array} \right.
$$
is not lower semicontinuous on $L^2(I)$.

In the formula above and throughout the paper we denote by $\gamma u$ the
trace of $u$ as a function from $BV(I)$ or $W^{1,p}(I)$, see \cite{Ziemer}.
Since $u\in BV(I)$ for   $I\subset \mathbb{R}$ may  have jumps, then by
definition  we have
$\gamma u(a) = \lim_{y\to a^+}u(y)$ and $\gamma u(b) = \lim_{y\to b^-}u(y)$.

\begin{definition}\label{d1}{\rm 
We shall say that a function $u\in L^2(0,T;L^2(I))$ is a solution to 
(\ref{rng}) if $u\in L^\infty(0,T;BV(I))$ and  $u_t\in L^2(0,T;L^2(I))$ and there is $\Omega \in L^2(0,T; W^{1,2})$. They satisfy the identity
\begin{equation}\label{rn-d1}
\langle u_t, \varphi \rangle = - \int_I \Omega \varphi_x\,dx
\end{equation}
for all test functions $\varphi \in C^\infty_0(I)$ and for almost every $t>0$.

\noindent
(A) We shall say that $u$ a solution to (\ref{rng}) satisfies the Neumann data
(\ref{rnN}), if
$$
\Omega|_{\partial I} =0\qquad\hbox{for }a.e.\ t>0.
$$
(B) We shall say that $u$ a solution to (\ref{rng}) satisfies the  Dirichlet
data (\ref{rnD}) at $x=a$ and $t>0$ if
$$
\gamma u(a) = A
$$
or
$$
\begin{array}{ll}
\hbox{if }  \gamma u(a) >A, &\hbox{then } \gamma\Omega(a) =2,\\
\hbox{if }  \gamma u(a) <A, &\hbox{then } \gamma\Omega(a) =-2.
\end{array}
$$
We shall say that $u$ a solution to (\ref{rng}) satisfies the  Dirichlet data
(\ref{rnD}) at $x=b$ and $t>0$ if
$$
\gamma u(b) = B
$$
or
$$
\begin{array}{ll}
\hbox{if }  \gamma u(b) >B, &\hbox{then } \gamma\Omega(b)=-2,\\
\hbox{if } \gamma u(b) <B, &\hbox{then } \gamma\Omega(b) =2.
\end{array}
$$}
\end{definition}
We notice that the time regularity postulated in Definition \ref{d1}
implies that solutions to  (\ref{rng}) are in $C([0,T];L^2(I))$. Hence,
we can impose initial conditions (\ref{rni}).

\noindent{\bf Remark.}
We also notice that our definition of solutions to (\ref{rng})  with Dirichlet boundary data coincides with that used by Andreu {\it et al.}, see \cite{dirichlet.andreu}.

We stress that $\Omega$ is a selection of the composition of the
multivalued operators $\cL\circ u_x$. 

We would like to expose the consequences of Definition \ref{d1}.
First, we will show the existence result. We note that we consider less regular initial conditions than in \cite{mury}.


\begin{theorem}\label{tw-gl}
Let us suppose that $u_0\in BV$, then\\
(1) there exists a unique solution to (\ref{rng}) with boundary conditions 
(\ref{rnD}), 
where $A,B\in \bR$;\\
(2) there exists a unique solution to (\ref{rng}) with boundary conditions 
(\ref{rnN});\\
(3) there exists a unique solution to (\ref{rng}) with periodic boundary
conditions (\ref{rnP}). \\
Moreover, 
for almost all $t>0$
\begin{equation}\label{rn-tw1}
\int_I [W(u_x + h_x) - W(u_x)]\,dx \ge \int_I \Omega  h_x\,dx,
\end{equation}
where $h\in C_0^{\infty}(I)$.
\end{theorem}

\noindent{\it Proof.} 
{\it Step 1.} After regularizing $\cL$ and $u_0$ we obtain a uniformly parabolic problem,
\begin{equation} \label{rnep}
 \begin{array}{ll}
  \frac{\partial u^\epsilon}{\partial t}= \left( {\cL^\epsilon}(u^\epsilon_x)\right) _x, & (x,t)\in I_T,\\
u^\epsilon(x,0)= u^\epsilon_0(x), & x\in I, \\
u^\epsilon(a,t)=A, \quad u^\epsilon(b,t)=B, & t>0,
 \end{array}
\end{equation}
where $\epsilon$ is a regularizing parameter. By the classical theory, see
\cite{lady}, we obtain existence and uniqueness of smooth solutions to
(\ref{rnep}). 

If we multiply (\ref{rnep}) by $u^\epsilon_t$ and integrate over $I_T$, then we reach,
$$
\int_0^T\int_I (u^\epsilon_t)^2\, dxdt = 
\int_0^T\int_I ( {\cL^\epsilon}(u^\epsilon_x))_x u^\epsilon_t\, dxdt.
$$
Integration by parts yields,
$$
\int_0^T\int_I (u^\epsilon_t)^2\, dxdt = - \int_0^T\int_I  {\cL^\epsilon}(u^\epsilon_x) u^\epsilon_{xt}\, dxdt = - \int_0^T\int_I \frac{d}{dt} W^\epsilon(u^\epsilon_x)\, dxdt,
$$
where $W^\epsilon(p)$ is the primitive of $\cL^\epsilon$ such that 
$$|p+1|+|p-1|\le W^\epsilon(p)\le |p+1|+|p-1|+k\epsilon$$
and $W^\epsilon(0)$ converges to  2 as $\epsilon\to 0$. 
Hence, we reach the following conclusion, 
\begin{equation} \label{rnee}
 \int_0^T\int_I (u^\epsilon_t)^2\, dxdt + \int_I W^\epsilon(u^\epsilon_x(x,T))\, dx =
\int_I W^\epsilon(u^\epsilon_{0,x})\, dx.
\end{equation}
Now, we will pass to the limit.
First of all, we notice that
$$\int_IW^{\epsilon}(u^{\epsilon}_{0,x}) \leq 3|I| + 2\sup_{\epsilon\in[0,1]}\int_I W(u^{\epsilon}_{0,x}) =: M.$$
Since we have found a bound on the right-hand-side (RHS) of (\ref{rnee})
independent of $\epsilon$ we conclude that
$$\int_0^T\int_I(u_t^\epsilon)^2\le M\quad \textrm{and}\quad \int_IW^\epsilon(u^\epsilon_x(x,t))\le M$$
for all $t\in [0,T]$. Thus, we can select a subsequence $\{u^\epsilon\}$ such
that
$$u^\epsilon\rightharpoonup u\quad \textrm{in}\quad L^2(0,T;L^2(I))$$
and 
$$u^\epsilon_t\rightharpoonup u_t\quad \textrm{in}\quad L^2(0,T;L^2(I)).$$
%
Furthermore, by  Aubin Lemma we deduce that
$$u^\epsilon\to u\quad \textrm{in}\quad L^p(0,T,L^q(I)),$$
where $p,q$ are arbitrary from the interval $(1,\infty)$. As a result for almost all $t\in(0,T)$
$$\|u^{\epsilon}(\cdot,t)-u(\cdot,t)\|_{L^q} \to 0.$$
Thus, we use the lower semicontinuity of the $BV$ norm to deduce that for almost all $t>0$
$$\varliminf_{\epsilon \to 0}\int_IW^{\epsilon}(u_{x}^{\epsilon})(x,t)dx \ge \varliminf_{\epsilon \to 0}\int_IW(u_{x}^{\epsilon})dx\ge \int_IW(Du)(\cdot,t),$$ 
where for $v\in BV (I)$ we write
\begin{equation} \label{rnW}
 \int_I W(D v)= \int_I | D(v+x)| + | D(v-x)|.
\end{equation}
Combining these inequalities, we arrive at
\begin{equation} \label{rne2}
 \int_0^t\int_I u _t^2(x,s)\, dxds + \int_I W(Du)(\cdot,t))\, dx \le M\quad \textrm{for almost all}\ t\in(0,T).
\end{equation}
We note that (\ref{rne2}) does not involve any statement on the boundary values
of $u$.

Moreover, we have a bound on the $BV$ norm of $u(\cdot,t)$,
\begin{eqnarray*}
 &\displaystyle{\int_I|Du|}&=\int_I|\frac{1}{2}D(u+x) + \frac{1}{2}D(u-x)|\\
& &\leq \frac{1}{2}\int_I|D(u+x)| + \frac{1}{2}\int_I|D(u-x)|=\frac{1}{2}\int_IW(u_x)dx.
\end{eqnarray*}

We also have to indicate a candidate for $\Omega$ as required by the definition
of a solution. We set 
$$
\Omega^\epsilon(x,t) := \cL^\epsilon(u^\epsilon(x,t)).
$$
Since $u_t^\epsilon = \Omega^\epsilon_x$, then due to (\ref{rnee}) we deduce that
\begin{equation}\label{rnom}
 \| \Omega^\epsilon \|_{L^2(0,T;H^1(I))} \le M_1<+\infty.
\end{equation}
Hence, we can select a  subsequence,
$$
\Omega^\epsilon\rightharpoonup \Omega\quad \textrm{in}\quad L^2(0,T;H^1(I)).
$$
Moreover,
$$
\int_0^T\int_I u_t \varphi\,dxdt = \int_0^T\int_I \Omega_x \varphi\,dxdt
\quad\hbox{for all } \varphi\in C_0^\infty((0,T)\times I).
$$
At this point we may apply \cite[Lemma 2.1]{mury} to conclude that
(\ref{rn-d1}) holds. Moreover, \cite[Lemma 2.2]{mury} implies (\ref{rn-tw1}). 

{\it Step 2.} We have to show that $u$ the limit of solutions to the regularized problems is a solution to 
(\ref{rng}) with boundary conditions, in the sense of Definition \ref{d1}. We first deal with more difficult case of Dirichlet boundary data (\ref{rnD}).

Let us suppose that $t$ is such that $u(\cdot, t)\in BV(I)$ and $\Omega(\cdot, t)\in W^{1,2}(I)$. We consider first $x=a$ a boundary point of $I$.
If $\gamma u(a,t) =A$, then $u$ satisfies the Dirichlet boundary  data.
Let us suppose that $A+\delta:=\gamma u(a,t) > A$. This means that for {\it any} sequence $\{x_n\}$, $a>x_n$, converging to $a$ we have $\lim_{n\to\infty}u(x_n,t) =A+\delta$. 
Then, we select $N$ such that for all $n>N$,
$$
A+ \frac12 \delta < u(x_n,t)=  u(x_n,t)- u^\epsilon (x_n,t) + u^\epsilon (x_n,t).
$$
Smooth solutions $u^\epsilon$ satisfy the boundary conditions in the above inequality, this implies that
$$
\frac12 \delta < u(x_n,t)- u^\epsilon (x_n,t) + u^\epsilon_x (c_n,t)(x_n-a),\qquad c_n\in(a,x_n).
$$
Since $u^\epsilon(\cdot,t)$ are commonly
bounded in $BV$, then we deduce with the help of Helly's theorem that
there is a subsequence  $u^\epsilon(\cdot,t)$, (we abstain from
introducing a new notation), that converges to $u(\cdot,t)$ everywhere.


We fix $n$ such that $\frac\delta{4(x_n-a)}>1$. Next, we 
select $\epsilon>0$ so that 
$u(x_n,t)- u^\epsilon (x_n,t)<\frac14 \delta.$ Combining these inequalities, we reach
$u^\epsilon_x (c_n,t)>1$, hence $\Omega^\epsilon(c_n,t) =2$. As a result, we conclude that
$$
\gamma\Omega(a,t)=2,
$$
as desired. The analysis of the remaining cases is similar, we leave it to the
interested reader. We conclude that $u$ is indeed a solution to (\ref{rng})
satisfying (\ref{rnD}).

The case of Neumann data is easier, since $\cL^\epsilon(u^\epsilon_x) \equiv
\Omega^\epsilon$ vanishes at ${\partial I}$ and the uniform convergence of
$\Omega^\epsilon$ as $\epsilon$ goes to zero implies that
$\gamma\Omega|_{\partial I}=0$. This happens for almost all $t>0$. 

The case of periodic boundary conditions (\ref{rnP}) is even easier than
(\ref{rnN}), because we do not have to worry about the values of $\Omega$.

{\it Step 3.} We shall establish uniqueness of solutions. The argument below is trivial in the case of  Neumann data. This is why we will consider only the Dirichlet boundary conditions. For the sake of simplicity of notation, we assume that  $[a,b] = [0,b]$.

We notice that if $u$ is a solution to (\ref{rng}), according to Definition \ref{d1}, then $t\mapsto u(t)\in L^2(I)$ is continuous, in particular it makes sense to evaluate $u$ at $t=0$.

Let us suppose that $u$ and $v$ are solutions to (\ref{rng}) satisfying
(\ref{rnD}) with $u(0)=u_0 = v(0)$. We will consider  functions $\tilde{u},
\tilde{v}$, which are extensions of $u,\ v$, respectively, to $\tilde{I} =
[-b,b]$. We set 
\begin{equation}\label{war-ok} 
\tilde{u}(x,t) = \left\{\begin{array}{cc}
                    u(x,t),& x\ge 0,\\
		    -u(-x,t), & x<0.
                   \end{array}\right.
\end{equation}
We define  $\tilde{v}$ in the same way. We also introduce extensions of
$\Omega(u)$ and $\Omega(v)$, where we stress the dependence of $\Omega$ 
on $u$ or $v$. We set,
$$
\tilde{\Omega} (x,t)= \left\{\begin{array}{cc}
                     \Omega(x,t),& x\ge 0,\\
		      \Omega(-x,t), & x>0.
                          \end{array}\right.
$$
We notice that  $\tilde{u}, \tilde{v}$,  $\tilde{\Omega}(\tilde u)$ and
$\tilde{\Omega}(\tilde v)$ may be extended as periodic functions with period
$2b$ to the line $\bR$. Thus, for any $\delta\in(0,2b)$, we have 
\begin{equation}\label{mnm}
\int_{\tilde{I}}\left(\tilde{\Omega}_x(\tilde{u}) -\tilde{\Omega}_x(\tilde{v})\right)(\tilde{u}-\tilde{v})\,dx = \int_{\tilde{I}+\delta}\left(\tilde{\Omega}_x(\tilde{u}) -\tilde{\Omega}_x(\tilde{v})\right)(\tilde{u}-\tilde{v})\,dx.
\end{equation}
Moreover, $\tilde{\Omega}(\tilde{u}) -\tilde{\Omega}(\tilde{v})$ belongs to $H^1(\tilde I)$, while
$\tilde{u}-\tilde{v}$ is a function from $BV(\tilde I)$. We want to perform
integration  by parts in the RHS of (\ref{mnm}). For this purpose, we quote the
following well-know lemma, see e.g. \cite{Loj}.
\begin{lemma}\label{cal-cze}
 If $w\in H^1(I)$, $\varphi\in BV(I)$, then
$$\int_Iw_x\varphi dx = -\int_ID\varphi + \gamma(w\varphi)|_a^b.\eqno\square$$
\end{lemma}
Due to the periodicity of the ingredients, we can select such $\delta$ that
$x=\delta$ is a point of continuity of
$\tilde{u}-\tilde{v}$. If we do so, we notice that the boundary terms drop out.
Thus, we conclude with the help of Lemma \ref{cal-cze} and formula (\ref{mnm})
that 
%
%
\begin{eqnarray*}
 \frac12 \|\tilde{u}-\tilde{v}\|_{L^2(\tilde{I})}^2(T) &= 
\displaystyle{-\int_0^T\int_{\tilde{I}+\delta} (\tilde{\Omega}(\tilde{u}; x,t)-\tilde{\Omega}(\tilde{v}; x,t))(\tilde{u}_x-\tilde{v}_x)\,dxdt}.
\end{eqnarray*}
On the other hand, monotonicity of $\cL$ yields
\begin{equation}\label{rn3}
 \frac12 \|\tilde{u}-\tilde{v}\|_{L^2(I)}^2(\tau) \le 0. 
\end{equation}
%
%
We conclude that $u=v$, as desired. \qed

\bigskip
The following observation is a consequence of the above proof and Lemma 2.1 in \cite{mury}. We make note of it for the future use.

\begin{corollary}\label{lm1}
(a) Let  us suppose that $u^n$, $\Omega^n $ is a sequence of  solutions to (\ref{rng}) such that $u^n\to u$ in $L^2(I_T)$ and
$\Omega^n \rightharpoonup \Omega$ in $L^2(0,T;H^1(I))$, then $u$ and $\Omega$ form a solution to (\ref{rng}).\\
(b) 
Let  us suppose that $u_0^n\in C^\infty$, $u_0^n\to u_0$ in $L^2(I) $ and $\sup \| u_0^n \|_{BV}<\infty$.
If $u_n \to u$ is a sequence of solutions to (\ref{rng}) with initial data $u_0^n$,
then,
$$
u_n \to u\quad\hbox{in } L^2(I_T),\qquad \Omega_n \rightharpoonup \Omega \quad\hbox{in } L^2(0,T;H^1(I).
$$ 
Hence, $u$ and $\Omega$ form a solution to (\ref{rng}) with initial data $u_0$.

\end{corollary}
{\it Proof.} Essentially, part (a) follows from the definition of a solution to (\ref{rng}). In order to show part (b) we recall that the a priori estimates we established in the course of proof of Theorem \ref{tw-gl} lead to a conclusion that 
$$
\| u^n \|_{L^\infty(0,T;BV(I))},\quad \| u^n_t \|_{L^2(I_T)},\quad \| \Omega^n \|_{L^2(0,T;H^1(I))}
$$
are uniformly bounded. Thus, we may extract a subsequence, not relabeled,  such that
$$
u_n \to u\quad\hbox{in } L^2(I_T),\quad \frac{\partial}{\partial t}u_n \rightharpoonup\frac{\partial}{\partial t}u\quad\hbox{in }  L^2(I_T),
\quad \Omega_n \rightharpoonup \Omega \quad\hbox{in } L^2(0,T;H^1(I).
$$
Due to part (a) $u$ is a solution to (\ref{rng}) with $u(\cdot,0) = u_0$. Since any from any subsequence we can select a converging subsequence to the same limit $u$, we conclude that the whole sequence $u^n$ converges to $u$, as well as $\Omega_n$ converges to $\Omega$. \qed

We also note a couple of estimates of solutions of (\ref{rng}), which will be used in next sections.

\begin{corollary}\label{col00}
 Suppose that $u_0\in BV$ and $u$ is the corresponding solution with either Dirichlet
(\ref{rnD}) for $A=B=0$ or Neumann (\ref{rnN}) data. Then, for all $t_2\ge t_1>0$, we have
$$
\| u(t_2) \|_{L^2} \le \| u(t_1) \|_{L^2} \le \| u_0 \|_{L^2}.
$$
\end{corollary}

\noindent{\it Proof.} We consider $u^\epsilon$ the solutions of regularized problem (\ref{rnep}). We test  this equation by $u^\epsilon$ and integrate over $I$ and $[0,t_1]$. We have
$$\int_0^{t_1}\int_Iu^\epsilon_tu\,dxdt = \int_0^{t_1}\int_I(\cL^\epsilon(u_x^\epsilon))_x u\,dxds.$$
We integrate by parts the RHS of the above equation, then
$$\int_0^{t_1}\frac{1}{2} \int_I\frac{d}{dt}(u^2)\,dxds = - \int_0^{t_1}\int_I(\cL^\epsilon(u^\epsilon_x))u^\epsilon_x\,dxds.$$
It follows from the monotonicity of $\cL$ that
$$\|u^\epsilon\|_{L^2(I)}(t_1)-\|u^\epsilon\|_{L^2(I)}(0)\le 0.$$
We obtain this inequality $\| u(t_2) \|_{L^2} \le \| u(t_1) \|_{L^2}$ by the replacing integration interval $[0,t_1]$ with $[t_1,t_2]$.
\qed

We can also prove in a similar way the following simple observation.
\begin{corollary}
 If $u$ and $v$ are two solutions to (\ref{rng}) satisfying Neumann, periodic or
the same Dirichlet
boundary conditions, then
\begin{equation}\label{rn-e}
 \|u-v\|_{L^2}(t) \le \|u_0-v_0\|_{L^2}. 
\end{equation}
\hfill \qed
\end{corollary}

The following estimate for $u_t$ is crucial for the rest of this paper.

\begin{proposition}\label{Prop1}
Suppose  that $u_0\in BV(I)$ and $u$ is the corresponding solution to
(\ref{rng})
with either (\ref{rnD}), (\ref{rnP}) or (\ref{rnN}) boundary conditions. Then, for a.e. $t\in (0, T)$, 
$$t\int_Iu_t^2(t,x)dx\leq \int_0^t\int_Iu_t^2(s,x)dxds\leq M<\infty.$$
\end{proposition}

\noindent{\it Proof.}
We proceed formally by differentiating equation (\ref{rng}) with respect to $t$ and testing it with
$u_t\varphi$, where $\varphi$ is non-negative and it depends only on $t$ and $\varphi(0)=0$. We have
$$u_{tt}u_t\varphi = \cL(u_x)_{xt}u_t\varphi.$$
Next, we integrate the above equation over $I$:
$$\frac{1}{2}\int_I\varphi\frac{d}{dt}u_t^2 = \int_I\cL(u_x)_{xt}u_t\varphi.$$
We integrate by parts the right hand side of the above equations, then
$$\frac{1}{2}\int_I\varphi\frac{d}{dt}u_t^2 = -\int_I\cL(u_x)_tu_{tx}\varphi.$$
Notice that, due to  monotonicity of $\cL$ we have $\cL(u_x)_tu_{xt}\varphi = \cL'(u_x)u_{xt}^2\varphi \geq 0$. Hence,
$$0\geq \int_I\varphi\frac{d}{dt}u_t^2dx.$$
We integrate the above equation over $[0,t]$, then
$$0\geq \int_0^t\int_I\varphi\frac{d}{dt}u_t^2\,dxds =
-\int_0^t\int_I\varphi'u_t^2dxds + \int_I\varphi u_t^2(t)\,dx$$
If $\varphi(t)=t$ and $u_0\in BV(I)$ then
\begin{equation}\label{I1}
t\int_Iu_t^2(t)dx\leq \int_0^t\int_Iu_t^2dxds\leq M<\infty.
\end{equation}
A rigorous argument is based on approximation. \qed
\bigskip

\begin{definition}{\rm 
We shall say that $t>0$ is {\it typical} if}
$$\int_I\left|{u_t(x,t)}\right|^2dx<\infty\quad \hbox{and}\quad \int_I\left|{\Omega_x(x,t)}\right|^2dx<\infty.$$
\end{definition}

\subsection{Discontinuous solutions}\label{sub2.1}

The type of convergence of approximate solutions permits discontinuous solutions. We make an observation about it.
 
\begin{proposition}\label{prdi}
Let us suppose that $u_0\in BV(I)$ and $u$ is the corresponding solution to  (\ref{rng}). If 
$u(\cdot, t_0)$ has a jump discontinuity at $x_0$ and $t_0$ is typical, then $\left|{\Omega(x_0,t_0)}\right| =2$.
\end{proposition}

\noindent{\it Proof.} Let us consider the  solutions $u^\epsilon$  to the regularized problem, approximating $u$. Since $u^\epsilon(\cdot,t)$ is a sequence of $BV$ functions, then by Helly Theorem we can select a subsequence (denoted by $u^\epsilon$) converging to $u$ everywhere. If $\Delta$ is the absolute value of the jump, then for a given 
$\epsilon$ we find $\delta_\epsilon$ such that
$$
\frac12 \Delta < |u^\epsilon(x_0+\delta_\epsilon,t_0) - u^\epsilon(x_0-\delta_\epsilon,t_0) |= 
2| u^\epsilon_x( c_\epsilon, t_0)| \delta_\epsilon .
$$
Thus, $| u^\epsilon_x( c_\epsilon, t_0)| >1$, as a result $|\Omega^\epsilon( c_\epsilon, t_0)|=2$. We can see that $|\Omega^\epsilon( c_\epsilon, t_0)| \to |\Omega(x_0,t_0)|=2$, because $ c_\epsilon \to x_0$. \qed

In the next section we will discuss steady states of (\ref{rng}). We will see that jump discontinuities of the solution are allowed also in steady states.

\section{Steady states}\label{set}

We describe the multitude of the steady states and we will consider all
three boundary 
conditions. We note that we frequently interchange $\Omega$ and $\cL(u_x)\in H^1(I)$. Here is our first observation.

\begin{proposition}\label{pr-st1}
(a) Let us suppose that a $BV$ 
function $u$ is a steady state solution to (\ref{rng}), i.e. there is $\Omega\in H^1$, understood as 
$\cL (u_x)$
satisfying
$$
(\cL (u_x))_x =0,
$$
then $\cL (u_x)$ is a constant from the set $\{ \pm 2, \pm 1, 0\}.$\\
(b) Let us suppose that  $u\in BV(I)$, but it is not Lipschitz continuous, and it satisfies 
$$
(\cL (u_x))_x =0,
$$
then $\cL (u_x)$ is a constant from the set $\{  2, -2\}$.
\end{proposition}

\noindent{\it Proof.} Part (a).\\
{\it Step i.} Of course, $\cL (u_x)$ is a constant from interval $[-2,2]$. Let
us suppose that $u_x>0$ on a set $E$ of a positive measure. Thus, on this set we
have $\sgn(u_x+1)=1$, as a result $\cL (u_x)\ge 0$ independently of the values
of $\sgn(u_x-1)$ on $E$. Let us suppose that $\cL (u_x)\in (1,2)$ on a set of
positive measure. We know that since  $\cL (u_x)$ is in $H^1$, then it is a
continuous function. Since
$$
1< \cL (u_x) = \sgn(u_x+1) + \sgn(u_x-1),
$$
and  $\sgn(u_x+1)=1$ on $E$, then $0< \sgn(u_x-1)<1$ and we  deduce that $u_x
=1$. Hence, any connected component of $E$ is a preimage of a facet, as a
result  $\sgn(u_x-1)$ may not be constant over $E$. 

Similarly we deal with the case $\cL (u_x)\in (-2,-1)$.

{\it Step ii.} Let us now suppose that $\cL (u_x)\in (0,1)$ on a set of positive measure. In this case we have
$$
0< \sgn(u_x+1) + \sgn(u_x-1) <1.
$$
But this implies an impossible situation of simultaneous $u_x+1=0$ and $u_x-1=0$ on the same set or
$\sgn(u_x+1)=1$ and $\sgn(u_x-1) <0$, i.e. $u_x+1>0$ and $u_x=1$. The last situation occurs on a facet, where $\sgn(u_x-1)$ may not be constant. 

Similarly, we deal with the case $\cL (u_x)\in (-1,0)$. 

Part (b) follows immediately from Proposition \ref{prdi} and part (a). 
\qed 

Our proposition states that the set of steady states may be very large. It should be stressed that it does not give a full description of this set, because the assumption is that $u\in BV$ conforms to Definition \ref{d1}.

\begin{proposition}\label{pr-st2}
 Let us suppose that a $u \in BV$ is a solution to (\ref{rng}) according to Definition  \ref{d1} and it is time independent, i.e., 
 function $u$ satisfies 
$$
(\cL (u_x))_x =0\qquad \hbox{in }(a,b).
$$
(a) If $\cL (u_x)= 0$ at $x=a$ and $x=b$, then $|u_x|\le 1$. That is, any Lipschitz function with the 
Lipschitz constant not exceeding 1 is a steady state of (\ref{rng}) with
homogeneous Neumann boundary conditions (\ref{rnN}).\\
(b) If $u$ is a steady state of (\ref{rng}) with  (\ref{rnD}) and $A\le B$ (the case $B\le A$ is analogous), then:\\
(i) if $(B-A)/(b-a) >1$, then any increasing function satisfying the data with
$u_x\ge 1$ is a steady state.\\
(ii) if $(B-A)/(b-a) =1$, then $u(x) = x+A-a$ is the only steady state.\\
(iii) if $(B-A)/(b-a) < 1$, then any Lipschitz continuous function with
$|u_x| \le 1$ satisfying $u(a)= A$, $u(b)=B$ is a steady state of (\ref{rng})
with (\ref{rnD}).\\
(c) If $u$ satisfies the periodic boundary condition, (\ref{rnP}), then $u$ is a
periodic Lipschitz continuous function with
$|u_x| \le 1$.
\end{proposition}

\noindent{\it Proof.} Part (a). 
Condition (\ref{rnN}) and Proposition \ref{pr-st1} jointly imply that
$\Omega(x) \equiv 0$. This, in turn yields that
$$
\sgn (u_x +1) = - \sgn(u_x -1)\neq 0.
$$
We note that the case $\sgn (u_x +1) = 0 = \sgn(u_x -1)$ is impossible. Thus,
$$
\sgn (u_x +1) =  1 = - \sgn(u_x -1).
$$
which implies
$$
u_x +1 \ge 0\quad\hbox{and}\quad u_x -1 \le 0
$$
or
$$
| u_x | \le 1.
$$
In particular, $u$ is Lipschitz continuous.

Part (b). Condition $A\le B$ implies that $u_x$ must be non-negative on a set
of positive measure. If $\cL (u_x)=2$, then $\sgn(u_x+1)=1 =\sgn(u_x-1)$, thus
$u_x\ge 1$. This implies that any monotone increasing function with $u_x \ge 1$
and 
$$
\lim_{x\to a^+} u(x) \ge A ,\qquad \lim_{x\to b^-} u(x) \le B
$$
is a steady state. In other words (i) holds.

If $\cL (u_x)=1$, then $\sgn(u_x+1)=1$ and $\sgn(u_x-1)=0$. This may occur only when $u_x=1$, i.e.
$u(x)= x+A-a$, thus $B=b-a+A$.

If $\cL (u_x)=0$, then $\sgn(u_x+1)=1$ and $\sgn(u_x-1)=-1$. This means that $u_x+1\ge 0$ a.e. and $u_x-1\le 0$ a.e. equivalently,
$$
|u_x| \le 1.
$$
In particular $u(a)=A$ and $u(b)=B$, for otherwise $\Omega$ would be equal to 2. 

Part (c). Of course, $\Omega$ may not be equal to $\pm 2$, because this would
imply that $u$ is increasing (or decreasing) on $I$, which is not possible for
a periodic  function. The argument presented above applies for the cases
$\Omega=\pm 1$, thus only $\Omega =0$ is left. As a result, we have the same
conclusion as in  the case of Neumann data.

\bigskip\noindent{\bf Remarks}\\
1) An important assumption is that $u$ is a solution to (\ref{rng}). We shall see that this imposes  restrictions of oscillations of $u_x$, see Theorem \ref{tmf}. \\
2) In particular, in case (iii) it seems  that an arbitrary number of
oscillations is possible.\\
3) If $(B-A)/(b-a)\le 1$, then all steady states are  Lipschitz continuous. On the
other hand, if  $(B-A)/(b-a) > 1$, then all increasing functions $u$, not
necessarily continuous, are steady states if $u_x\ge 1$ (see Proposition
\ref{prdi}). Thus, we see that for $u_0\in BV(I)$ the regularity $u\in
L^\infty(0,T;BV(I))$ is optimal.

\section{Properties of solutions}\label{sep}

We collect properties of solutions related to facets, defined in \S \ref{serg}, and their evolution. In \S \ref{sect} we study extinction times. 

\subsection{Facets}\label{serg}
In this section we will answer the question if solutions may
have infinitely many facets. We will introduce necessary notions.

\begin{definition} \label{2}{\rm Let us set $\cP = \{-1, 1\}$.\\
 (a) We shall say that a subset $F$ of the graph of a solution
to (\ref{rng}) is a {\it facet}, if 
$$
F \equiv F(\xi^-,\xi^+)=\{ (x,u(x)):\ u_x|_{[\xi^-,\xi^+]} \equiv p \in \cP
\}.
$$
We write $u_x|_{[\xi^-,\xi^+]}$ with the understanding that the one-sided
derivatives of $u$ exist at $\xi^+$ and $\xi^-$.

Moreover, if $[\xi^-,\xi^+]\subset J$, $J$ is an interval  and $u_x|_J\equiv p\in \cP$, then $[\xi^-,\xi^+] =J$. Interval  $[\xi^-,\xi^+]$ is called the
pre-image of facet $F$ or a faceted region of $u$, (cf. Section \ref{lepkie}).

(b) Facet $F(\xi^-,\xi^+)$ has {\it zero curvature}, if: (i)
either $\Omega$ i.e. $\cL(u_x)$ have the same value at $\xi^-$ and $\xi^+$ or (ii) $\xi^- =a$ or $\xi^+ =b$ i.e. the facet hits the boundary (in the case of Dirichlet boundary conditions).}
\end{definition}

\begin{remark}
 We notice that if for a facet $F=F(\xi^-,\xi^+)$ there is $\delta$ such that $u\left|_{(\xi^--\delta,\xi^+ +\delta)}\right.$ is monotone, then $F$ has zero curvature.
\end{remark}
This will be seen in the proof of Lemma \ref{lem-chi}.

We begin our analysis with the following observation.
\begin{lemma}\label{leb}{\rm 
Let us suppose that  $t>0$ is typical 
and $F(\xi_l,\xi_r)$ is a facet, then 
 $$
 \lim_{x\to\xi_l^-} \Omega(x,t) = \omega^-\in\{-2,0,2\},\qquad
 \lim_{x\to\xi_r^+} \Omega(x,t) = \omega^+\in\{-2,0,2\}.
 $$
}\end{lemma}

\noindent{\it Proof.} We shall investigate the neighborhood of $\xi_r$ assuming that $ u(x,t)$ is absolutely continuous. Then, there are two possibilities, (1) for all sufficiently $\epsilon>0$, 
the derivative $u_x$ on $(\xi_r, \xi_r+\epsilon)$ assumes  values from the set $\{-1,1\}$; (2) there is a sequence $\{ x_n\}_1^\infty$ converging to $\xi_r$ such that $\xi_r<x_n$ and $u_x(x_n,t)$ exists and $u_x(x_n,t)\neq \pm 1$.

If the first case occurs and there is $\delta>0$ such that for all $x\in (\xi_r, \xi_r+\delta)$ we have 
$u_x(x,t) = - u_x|_{[\xi_l,\xi_r]}$, then one can see that
$$
\lim_{x\to \xi^+_r} \Omega(x,t) =0.
$$

If instead there are two sequences $x_n^+,x_n^-$ converging to $\xi_r$, such that $\xi_r<x_n^+,x_n^-$  and $u_x(x_n^+,t) =1$,  $u_x(x_n^-,t) = -1$, then
$$
\sgn(u_x(x_n^+,t)+1) + \sgn(u_x(x_n^+,t)-1) = 1 +\zeta^+_n,$$
$$ 
\sgn(u_x(x_n^-,t)+1) + \sgn(u_x(x_n^-,t)-1) = \zeta^-_n - 1.
$$
Since the right limit of $\Omega$ exists at $\xi_r$ we deduce that
$\zeta^\pm := \lim_{n\to \infty} \zeta^\pm_n$ satisfy
$$
2 = \zeta^- - \zeta^+.
$$
As a result, $ \zeta^- = 1 =- \zeta^+$ and 
we conclude that $\Omega(x_n^+,t) = 0 = \Omega(x_n^-,t)$ and 
$$
\lim_{x\to \xi^+_r} \Omega(x,t) =0.
$$

Now, we  consider the situation, when there is a sequence $x_n$ converging to $\xi_r$, such that $\xi_r<x_n$ and $u_x(x_n, t)$ exists and $u_x(x_n, t)\neq \pm 1$. If this happens, then (a) $u_x(x_n ,t) > 1$, (b)  
$u_x(x_n ,t) \in (-1,1)$ or (c) $u_x(x_n ,t) < - 1$. Case (a) leads to the conclusion, that 
$\Omega(x_n ,t) = 2,$ hence 
$$
\lim_{x\to \xi^+_r} \Omega(x,t) =2.
$$
If (b) takes place, then $\Omega(x_n ,t) = 0$  and
$$
\lim_{x\to \xi^+_r} \Omega(x,t) =0.
$$
Finally, (c) leads to $\Omega(x_n ,t) = -2$. As a result,
$$
\lim_{x\to \xi^+_r} \Omega(x,t) =-2.
$$

If we consider general data, i.e. $u(\cdot,t)\in BV(I)$ for a.e. $t>0,$ then this solution may be approximated by smooth solutions $u^n$ for which the above conclusion is valid. We notice that sequence $\Omega^n$ converges uniformly, hence our claim follows for $\xi_r$. A similar conclusion may be drawn for $\xi_l$. \qed

\bigskip
We immediately deduce that: 
(see also \cite{mury})
\begin{proposition}\label{pr-b}
For a typical $t>0$ $u(t)$ does  not contain any degenerate $(\xi^-=\xi^+=:\xi)$ facet with nonzero  curvature.
\end{proposition}

\noindent{\it Proof.} 
Let us suppose the contrary, i.e. there is  $\xi\in I$ such that $u_x(t,\xi)=1$.
Then, according to Lemma \ref{leb}, $\lim_{x_n\to \xi^\pm}=\Omega(x_n,t) \in \{-2,0,2\}$. However, the continuity of $\Omega$ implies that
$$\Omega(\xi^-,t) = \Omega(\xi^+,t).$$
Hence $F(\xi)$ has zero curvature. \qed

We present our main structural theorem.
\begin{theorem}\label{tmf}
 If $u_0\in BV(I)$ and $u$ is the corresponding solution to (\ref{rng}) with
either Dirichlet (\ref{rnD}), periodic (\ref{rnP}) or Neumann (\ref{rnN}) boundary conditions, then
for almost all $t>0$ the number of facets with curvature different from zero is
finite. 
\end{theorem}

We will prove the Theorem for a typical $t>0$. 
The proof is based on the following estimate for facet velocity.

\begin{lemma}\label{lm2}
If $u$ is a solution to (\ref{rng}) and $F(\xi^-,\xi^+)$ is one of the facets,
then 
$\Omega_x=\sgn (u_x-1)_x +\sgn (u_x + 1)_x$ has the following form,
\begin{equation}\label{rn-4.1}
 \Omega_x =  \frac{\omega^+- \omega^-}{\xi^+-\xi^-},
\end{equation}
where $\omega^\pm:= \Omega(\xi^\pm,t)$.
\end{lemma}

\noindent{\it Proof.} {\it Step 1.}\\
Let us denote by $E$ one of the functionals on $L^2$
$$
E_1(u)=\left\{
\begin{array}{ll}
 \int_I W(Du) & u\in BV(I),\\
+ \infty & u \in L^2(I)\setminus  BV(I),
\end{array}
\right. 
E_2(u)=\left\{
\begin{array}{ll}
 \int_{\bT} W(Du) & u\in BV(\bT),\\
+ \infty & u \in L^2(\bT)\setminus  BV(\bT),
\end{array}
\right.
$$
where $\bT$ is a flat one-dimensional torus identified with $[0, b)$. Functional $E$ is convex
and lower semicontinuous on $L^2$. Thus, we may consider the gradient flow
\begin{equation}\label{rngl}
 u_t \in - \partial E(u),\quad u(0)= u_0.
\end{equation}
By the theory of nonlinear semigroups, see \cite{brezis}, we know that for any
$u_0\in L^2$, there exists $v$, a unique solution to (\ref{rngl}), such that
$v(t) \in D(\partial E(v(t)))$, $\frac {dv^+}{dt} $ exists for all $t>0$ and
$$
\frac {dv^+}{dt} = \Upsilon_x,
$$
where $\Upsilon_x$ is the minimal section of $\partial E(v)$.

It is easy to check that $\Omega(\cdot, t)\in \partial E(u(t))$, (see (\ref{rn-tw1}), where $u$ is
constructed in Theorem \ref{tw-gl}, see also \cite{mury}.
Moreover, 
$$t \to \int_I W(v_x(x,t))$$
is a decreasing function on $[\delta,\infty)$ for $\delta >0$, provided that $v_0\in L^2(I)$ (see \cite{brezis}). In particular, $v(\cdot,t)\in BV(I)$ for almost all $t>0$. 
Thus, (see Lemma \ref{cal-cze}),
$$
\frac d{dt} \int_I |u-v|^2 \,dx = \int_I (\Omega_x- \Upsilon_x)(u-v) \,dx 
= -\int_I (\Omega- \Upsilon)(u-v)_x\,dx\le 0.
$$
We notice that due to the Neumann or periodic boundary conditions the boundary term vanishes. Hence, 
$$
u(t) = v(t)
$$
for all $t\ge 0$. Now, we can use the fact that $\Omega$ is the minimal section
of $\partial E(u)$. Thus, if $F(\xi^-,\xi^+)$ is a facet, then 
$\Omega$ minimizes the functional
$$
\int_{\xi^-}^{\xi^+} | \zeta_x|^2\,dx
$$
among $H^1$ functions
with the specific boundary conditions, which we will discover momentarily. As a result,
$\Omega$ is a linear function. From Lemma \ref{leb} we know the values of $\Omega$ at the endpoints of facet $F(\xi^-,\xi^+)$. In particular, $\Omega(\xi^-,t) \neq \Omega(\xi^+,t) $ if and only if $F(\xi^-,\xi^+)$ is not of zero curvature. Thus (\ref{rn-4.1}) holds.

{\it Step 2.} Now,  we study eq. (\ref{rng}) with Dirichlet data and let us suppose that $u$ is a
solution to this problem. For the sake of simplicity we may assume that
$I=[0,b]$, then we may extend $u$ to an odd function $\tilde u$ on $[-b,b]$. Now, it is convenient to
identify $[-b,b]$ with $\bT$. Thus,  we may apply results of Step 1. 
In particular, if $F(\xi^-,\xi^+)$ is a facet, then $\Omega$ restricted to $[\xi^-,\xi^+]$ is an affine function.
\qed

\bigskip
We return to the {\it proof of Theorem \ref{tmf}}.
By Proposition \ref{Prop1} we know that for
almost all $t>0$ we have
$$
\int_I u_t^2\,dx \le \frac{1}{t}\int_0^t\int_I u_t^2\,dxdt <\infty.
$$
Thus, we square the RHS of (\ref{rng}) and integrate $
u_t^2= |\Omega_x|^2. 
$ over $I$.
We notice that
\begin{eqnarray*}
\int_I u_t^2\,dx = \int_I|\Omega_x|^2\,dx \ge 
 \sum_{F(I_\iota)} \int_{I_\iota} |\Omega_x|^2\,dx 
 =\sum_{F(I_\iota)} \frac{(\omega^+-\omega^-)^2}{\xi^+_i-\xi^-_i}.
\end{eqnarray*}
Here, $\{F(I_\iota):\ \iota \in J\}$ is the collection of all non-zero cuvature facets.
We immediately conclude that the number of facets with non-zero curvature is
finite. \qed

%
%
%


\noindent{\bf Remark.} The extra effort for the Dirichlet boundary data is related to the fact that functional $E_1$ is not lower semicontinuous in $L^2$ with the natural choice of domain, i.e. 
$\{u\in BV:\ \gamma u(a) = 0 = \gamma u(b)\}$.

After these preparations we are going to present the basic facts about facet creation process. Our main tool is analysis of continuity of $\Omega(\cdot,t)$. Let us first notice.

\begin{lemma}\label{lemA}
 If $u$ is a solution to (\ref{rng}), $t>0$ is a typical time instance and $F=F(\xi^-,\xi^+)$ is a non-degenerate facet, which does not touch the boundary i.e. $\xi^-\neq a$ and $\xi^+\neq b$, then $|\Omega(\xi^+,t) - \Omega(\xi^-,t)| \le 2$.
\end{lemma}

\noindent{\it Proof.} In order to fix our attention we consider $u_x = p =1$ on $(\xi^-,\xi^+)$. Then, for any 
$x\in (\xi^-,\xi^+)$ we have
$$
\Omega(x,t) = \sgn( u_x(x,t) +1) + \sgn( u_x(x,t) -1) = 1 + \zeta(x,t),
$$
where $\zeta(x,t) \in [-1,1]$. Hence, 
$$
|\Omega(\xi^+,t) - \Omega(\xi^-,t)| =  |\zeta(\xi^+,t) - \zeta(\xi^-,t)|
\le 2.
$$
Similar analysis is valid for $p=-1$. Our claim follows. \qed


The following fact is quite important, it tells us that a solution $u$ does not miss any of the preferred directions $-1, 1$, even if the datum does.

\begin{theorem}\label{proB}
 Let us suppose that  $u$ is a solution to (\ref{rng}), $t>0$ is a typical time instance,
i.e. $u_t(\cdot, t), \Omega_x(\cdot, t) \in L^2(I)$. If $x_0\in I$ and $u_x^-(x_0, t)< u_x^+(x_0, t)$, (resp.  $u_x^+(x_0, t)< u_x^-(x_0, t)$), then
$$
(u_x^-(x_0, t), u_x^+(x_0, t)) \cap \{-1, +1\} = \emptyset, \qquad(\hbox{resp. }
(u_x^+(x_0, t), u_x^-(x_0, t)) \cap \{-1, +1\} = \emptyset).
$$
\end{theorem}

\noindent{\it Proof.} First, we consider the case of $u(\cdot, t)$ being absolutely continuous.

 Let us assume that the opposite happens, i.e. there is $p$ from $\{-1,1\}$ such that
$$p\in (u_x^-(x_0,t), u_x^+(x_0,t)).$$
For the sake of definiteness we assume that $p=1$. In other words,
$$
u_x^-(x_0,t) <1 < u^+_x(x_0,t).
$$
Thus, for all $x>x_0$ sufficiently close to $x_0$ we have
$$
1< \frac{u(x,t) -u(x_0,t)}{x-x_0}.
$$
We conclude that there exists a sequence $x^+_n$, converging to $x_0$ such that $x_0<x^+_n$ and $1<u_x(x^+_n,t)$. Hence,
$$
\Omega(x^+_n,t) = \sgn( u_x (x^+_n,t) + 1) +  \sgn( u_x (x^+_n,t) -1) = 2.
$$
Continuity of $\Omega$ at $x_0$ implies that
$$
\Omega^+(x_0,t) := \lim_{x_n^+\to x_0^+} \Omega(x_n^+,t) =2
$$
On the other hand, $u_x^-(x_0,t) <1$ and we see that
$$\frac{u(x,t) -u(x_0,t)}{x-x_0}<1\quad \textrm{i.e.}$$
there exists a sequence $x^-_n<x_0$ converging to $x_0$ such that $u_x(x^-_n,t)<1$. As a result
$$\Omega(x_n^-,t) = \sgn( u_x (x^-_n,t) + 1) +  \sgn( u_x (x^-_n,t) -1) = \zeta_n-1.
$$
Now, we consider three cases depending on the behavior of $u_x(x_n^-,t)$. We
set $\Omega^-(x_0,t)=\lim_{x_n^-\to x_0}\Omega(x_n^-,t)$.
Now if $u_x(x^-_n)<-1$, then $\zeta_n =-1$ and $\Omega^-(x_0,t) = -2$,
if $u_x(x^-_n)>-1$, then $\zeta_n=1$ and $\Omega^-(x_0,t) =0.$ If $u_x(x^-_n)=-1$, then $\zeta_n\leq 1$ and then $\Omega^-(x_0,t)\leq 0$.

Let us consider the possibility that
$u_x^+=1$, then $u_x^-<1$, but $\frac{u(x)-u(x_0)}{x-x_0}>1$, thus
$\Omega^+(x_0,t) = 1+\zeta^+=2$ and $\Omega^-(x_0,t) = -1+\zeta^-\leq 0$.

Let us consider the case of general data, i.e. $u(\cdot,t)\in BV(I)$. We may assume that $u(\cdot,t)$ has jump discontinuity at $x_0$. It follows from Proposition \ref{prdi} that $\left|{\Omega(x_0,t)}\right|=2$.  In order to fix attention we assume that $\Omega(x_0,t)=2$. By continuity of $\Omega$, there is a neighborhood $U$ of $x_0$ such that $\Omega(x,t)>2-\varepsilon$ for $x\in U$. Hence, there does not exist any $x\in U$ such that $u_x(x,t)<1$. Therefore, $u_x(x,t)\geq 1$ for all $x\in U$. Thus, the open interval with endpoints $u_x^-, u_x^+$ does not contain $+1$ nor $-1$. A similar result holds for $\Omega(x_0,t)=-2$.   
\qed

\bigskip
We introduce a piece of convenient notation. Let us suppose that $t>0$ is typical and $u(\cdot,t)\in BV(J)$, $J=[\alpha,\beta]$, has a facet $F(\xi^-,\xi^+)$, we assume that $\alpha<\xi^-\le \xi^+<\beta$. We introduce the {\it transition  numbers}  $\chi_s= \chi_s(u ,x)$, $s=l,r$, by the following formula,
\begin{equation}\label{tr-n1}
 \chi_{l}=\begin{cases}
         +1 \ \ {\rm if} \ u \geq \ell_{p} \ \ {\rm in} & \{ x\in J : x\leq
x_{0} \}, \\
         -1 \ \ {\rm if} \ u \leq \ell_{p} \ \ {\rm in} & \{ x\in J : x \leq
x_{0} \}, 
\end{cases}
\end{equation}

\begin{equation}\label{tr-n2}
\chi_{r}=\begin{cases}
         +1 \ \ {\rm if} \ u \geq \ell_{p} \ \ {\rm in} & \{ x\in J : x\geq
x_{0} \}, \\
         -1 \ \ {\rm if} \ u \leq \ell_{p} \ \ {\rm in} & \{ x\in J : x
\geq x_{0} \},
\end{cases}
\end{equation}
where $ \ell_{p}$ is the line with slope $p$ containing  facet $F(\xi^-,\xi^+)$.

We note that if facet $F(\xi^-,\xi^+)$ has zero curvature, then $\chi_l(u ,x)+\chi_r(u ,x) = 0.$

\begin{lemma}\label{lem-chi}
 If $u$ is a solution to (\ref{rng}) and $F(\xi^-,\xi^+)$ is one of the facets, $\xi^-\neq a$ and  $\xi^+-\neq b$, 
then for a typical $t>0$,
\begin{equation}\label{omega-chi}
 \Omega_x = \frac{\chi_l+\chi_r}{\xi^+-\xi^-}.
\end{equation}

\end{lemma}

\noindent{\it Proof.}
We are going to find values of $\Omega$ at $\xi^-$, $\xi^+$. We treat $\xi^+$ first.
Initially, we assume that $u(\cdot,t)\in AC(I)$ and $t>0$ is a typical time instance.
Moreover, at $t>0$ 
$u(\cdot,t)$ does not have any degenerate facets, as guaranteed by Proposition \ref{pr-b}. Thus, for facet $F(\xi^-,\xi^+)$ there is such $\epsilon>0$ that $u|_{(\xi^+, \xi^+ +\epsilon)}$ is either above $l_p$, i.e. the line containing  $F(\xi^-,\xi^+)$, or below it.

If $u$ is above $l_p$, then 
$$
u(x,t) - u(\xi^+,t) > p(x -\xi^+)\qquad \hbox{for all } x\in (\xi^+, \xi^++\epsilon), 
$$
where $p$ is the slope of $l_p$. This implies that there exists a sequence $x_n \in (\xi^+, \xi^+ +\epsilon)$, converging to $ \xi^+$ such that $u_x(x_n)>p$. We notice that if $p=1$, then
$$
\Omega(x_n,t) = \sgn(u_x(x_n,t)+1) + \sgn(u_x(x_n,t)-1) = 1 + 1 = p + \chi_r.
$$
If  $p=-1$, then we know from Theorem \ref{proB} that $u_x(x_n,t)\in (-1,1)$. Thus we see that
$$
\Omega(x_n,t) = \sgn(u_x(x_n,t)+1) + \sgn(u_x(x_n,t)-1) = 1 - 1 = p + \chi_r.
$$
Since $\Omega(\cdot, t)$ is continuous, we conclude that
\begin{equation}\label{rn-chip}
 \Omega(\xi^+,t) =  p + \chi_r.
\end{equation}
A similar reasoning performed for interval $(\xi^- -\epsilon, \xi^-)$ yields
\begin{equation}\label{rn-chil}
 \Omega(\xi^-,t) =  p - \chi_r,
\end{equation}
provided that $u(\cdot, t)$ is continuous. 

Let us suppose now that $u(\cdot, t)$ is no longer absolutely continuous in any neighborhood of $x_0$. What may happen is:\\
(a) $u(\cdot, t)$ has a jump discontinuity at $x_0$, thus by Proposition \ref{prdi} we know that $\Omega(x_0,t) = 2$.\\
(b)  $u(\cdot, t)$ is  continuous at $x_0$. We consider 
$$
u(x_n) - u(x_0) = p(1+\delta_n)  (x_n-x_0)> p (x_n-x_0).
$$
If $u^\ep$ is the regularized solution, then we can find a sequence $\epsilon_n$ converging to 0 such that
$$
u^{\epsilon_n}(x_n) - u^{\epsilon_n}(x_0) \ge  p(1+\frac12\delta_n)  (x_n-x_0)> p (x_n-x_0).
$$
Thus, we may use the argument from the  first part. 
However, even if $u^k(\cdot, t)$ converges to  $u(\cdot, t)$ in $L^2$, then $\Omega^k(\cdot,t) $ converges uniformly to  $\Omega(\cdot,t)$. Thus, (\ref{rn-chip}) and (\ref{rn-chil}) remain valid for solutions with $u_0\in BV$. Thus, in all cases, mentioned above, we obtain that
\begin{equation}\label{rn-om}
 \Omega_x(x,t) = \frac{\Omega(\xi^+,t) -\Omega(\xi^-,t) }{\xi^+ - \xi^-}
\end{equation}
for $x\in (\xi^-,\xi^+)$, i.e. (\ref{rn-4.1}) holds.

\qed

The two previous results are concerned with typical time instances. In particular, they do not preclude the possibility of shrinking a non-zero curvature facet to a point at an exceptional time. Now, 
we present an improvement of Theorem \ref{tmf}. The Proposition below is
not a direct consequence of regularity.

We used $\xi^-$, $\xi^+$ to denote the endpoints of the facet pre-image. Now, it is advantageous to show the dependence of $\xi^-,\xi^+$ on $h$ the \textquoteleft distance\textquoteright of the facet from the $x_1-$ axis. If $t_0>0$ is a typical time instance, then we set $u(\cdot)=u(\cdot,t_0)$ and $h_0=u(x_0,t_0)$. We have
%
%
\begin{equation}\label{defXi-}
 \xi^- = \inf\{x:u(x) = p(x-x_0)+h\}, 
\end{equation}
\begin{equation}\label{defXi+}
 \xi^+=\sup\{x:u(x) = p(x-x_0)+h\}.
\end{equation}
We notice that $\xi^\pm$ are well-defined. Moreover, due to the argument in the above proof,
\begin{equation}\label{warXI}
 \begin{array}{ccc}
  \frac{d\xi^+}{dh} \geq 0\quad \textrm{and}\quad \frac{d\xi^-}{dh} \leq 0 & \textrm{if}& \chi_l+\chi_r >0,\\
\frac{d\xi^+}{dh} \leq 0\quad \textrm{and}\quad \frac{d\xi^-}{dh} \geq 0 & \textrm{if}& \chi_l+\chi_r <0.
 \end{array}
\end{equation}

\begin{proposition}\label{cc}
 No facet with non-zero curvature may shrink to a point at any $t>0$.
\end{proposition}

\noindent{\it Proof, step 1.} We begin with one facet $F(\xi^-,\xi^+)$ with non-zero curvature. We assume that $F(\xi^-,\xi^+)$ does not intersect other non-zero curvature facets. Let
$$L=\xi^+-\xi^-.$$
After taking the time derivative, we have
$$\frac{d}{dt}(\xi^+-\xi^-) = \left(\frac{d\xi^+}{dh}-\frac{d\xi^-}{dh}\right)\frac{dh}{dt}.$$
We notice that due to (\ref{warXI}) the RHS is always non-negative.

{\it Step 2.} Now we consider two  intersecting  facets  $F(\xi^-,\xi)=F_1,$ $F(\xi, \xi^+)=F_2$ with non-zero curvature (see Figure 1), i.e $(\xi,u(\xi))$ is the intersection point. 
\begin{figure}[h]
 \begin{center}
{\epsfxsize=8cm \epsfbox{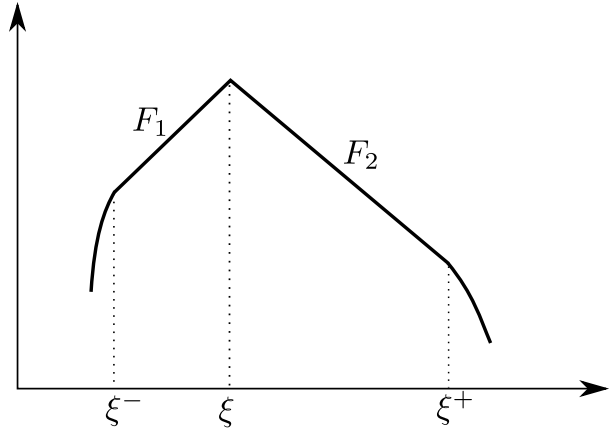}}
\end{center}
\caption{Two interacting facets.}
\end{figure}

For facet $F_1$ it is advantageous to use $\xi^-$ defined by (\ref{defXi-}) and for $F_2$ to use $\xi^+$ defined by (\ref{defXi+}).
Point $\xi$ is the intersection of lines containing $F_1$ and $F_2$, i.e.,
$$\xi = \frac{x_1+x_2}{2}+ \frac{h_2-h_1}{2p},$$
where $(x_i,h_i)\in F_i$, $i=1,2$ are fixed and $p$ is the slope of $F_1$.
If we set $L_1=\xi - \xi^-$, $L_2=\xi^+-\xi$, then we see that
$$\frac{d}{dt}{L_1} = \frac{d\xi}{dt} - \frac{d\xi^-}{dh_1}\frac{dh_1}{dt}=
\frac{1}{2p}\left(\frac{dh_2}{dt}-\frac{dh_1}{dt}\right)- \frac{d\xi^-}{dh_1}\frac{dh_1}{dt},$$
$$\frac{d}{dt}{L_2} =
\frac{1}{2p}\left(\frac{dh_1}{dt}-\frac{dh_2}{dt}\right)+ \frac{d\xi^+}{dh_2}\frac{dh_2}{dt}.$$
We notice that the last terms are always non-negative.
Moreover from the formula for vertical velocity, (\ref{omega-chi}), we calculate $\frac{dh_1}{dt}$, $\frac{dh_1}{dt}$ and we have
$$\frac{d}{dt}{L_1} = \frac{\chi_l+\chi_r}{2p}\left(\frac{L_1-L_2}{L_1L_2} \right)- \frac{d\xi^-}{dh_1}\frac{dh_1}{dt},$$
$$\frac{d}{dt}{L_2} = \frac{\chi_l+\chi_r}{2p}\left(\frac{L_2-L_1}{L_1L_2} \right)+\frac{d\xi^+}{dh_2}\frac{dh_2}{dt}.$$

We notice that if $L_1\to 0$ and $L_2\geq \delta >0$, then  we see that $\frac{dL_1}{dt} >0$,
also  if $L_2\to 0$ and $L_1\geq \delta >0$, then $\frac{dL_2}{dt} >0$, but this 
is impossible. 
Finally we notice that $L_1+L_2\to 0$ is impossible too, because
$$\frac{d}{dt}(L_1+L_2) = \frac{d\xi^+}{dt} -\frac{d\xi^-}{dt}\geq 0.\eqno\Box$$

We see that only zero curvature facets may shrink to a point. On the other hand zero curvature facets are created during  collisions.


\subsection{Extinction time}\label{sect}
We notice that the diffusion is so strong that for all initial data $u_0$ the
solution gets extinct in finite time. The examples from Section \ref{ssdisc} give
explicit  bounds in the case of initial data $u_0\in BV(I)$. 

\begin{definition}{\rm 
 If $u$ is a solution to (\ref{rng}), then a number $T_{ext}>0$ is the
extinction time for $u$ iff $u_t \equiv 0$ for all $t>T_{ext}$ and for all
$\epsilon>0$ we have $u_t \not\equiv 0$ on $(T_{ext}-\epsilon,T_{ext}).$}
\end{definition}


%


First we establish the following proposition, which will be helpful in the
next theorem.

\begin{proposition}\label{prop-Text}
Suppose that 
$v$ and  $u^k$, $k\in \mathbb{N}$ are solutions to (\ref{rng}),  $T_{ext}^k$  is the extinction time of $u^k$. If $u^k \rightarrow v$
in $L^2(I_T)$,  $\Omega^k \rightharpoonup \Omega$ and $ T_{ext}$ is the extinction time of  $v$,
then
$$
T_{ext} \le \limsup_{k\rightarrow\infty} T_{ext}^k.
$$
\end{proposition}

\noindent{\textit{Proof}.} Set $\bar T := \limsup_{k\rightarrow\infty} T_{ext}^k$. It is sufficient to show that if $\delta>0$ and $h\in\bR$ are such that $|h| <  \delta$, then
$$
v( \bar T  +\delta + h) - v( \bar T  +\delta ) =0.
$$
Indeed, since $u^k\in C([0, \bar T  +\delta + h]; L^2(I))$ and $u^k$ converges uniformly to $v$ in  
$C([0, \bar T  +\delta + h]; L^2(I))$, then we have
\begin{equation}\label{rn4}
 v( \bar T  +\delta + h) - v( \bar T  +\delta ) =
\lim_{k\rightarrow \infty}\left[ u^k ( \bar T  +\delta + h) - u^k( \bar T  +\delta )\right] =0.
\end{equation}
The limit is in $L^2(I)$. 
This is so, because for sufficiently large $k\in \mathbb{N}$, we have
$T_{ext}^k < \min\{ \bar T  +\delta + h, \bar T  +\delta\}$. This implies that
$u^k ( \bar T  +\delta + h)= u^k( \bar T  +\delta ) =  u^k( T_{ext}^k)$. Hence, (\ref{rn4}) follows. \qed

\begin{theorem}\label{T-ext} Let us suppose that $\delta>0$ is any typical time. We set
$u_0:= u(\delta)$, where $u$ is a solution to (\ref{rng}). We assume that $u_0$ is differentiable away form the endpoint of the facets.  If $x_1$ and $x_2$ are points with the same value of the derivative of $u_0$,
either $+1$ or $-1$ and there is no $x_0\in[x_1,x_2]$ with the opposite value
of derivative of $u_0$, then $T(x_1,x_2)$ the time after which the facets at $x_1$
and $x_2$ collide is finite. Then the extinction time $T_{ext}$ of $u_0$
can be estimated as 
$$\begin{array}{cl}
T_{ext}
 \le \max\{ T(x_1,x_2):& (u_0)_x(x_1)=(u_0)_x(x_2)\\
&\textrm{there is no } x_0\in[x_1,x_2]\textrm{ such that }  (u_0)_x(x_0) = -(u_0)_x(x_1)\}.
\end{array}
$$
\end{theorem}

\noindent{\it Proof.} By assuming that $\delta>0$ is typical, 
in virtue of Proposition \ref{pr-b} and Proposition \ref{cc} we consider situation, when all facets have been already created and their lengths are positive. It follows from Theorem \ref{tmf} that for $t>\delta>0$ we have finitely many facets with non-zero curvature. 
Let
$$\mathcal{F}^+ = \{F(\xi^-,\xi^+):u_x|_{[\xi^-,\xi^+]}=1, 
\chi_l(\xi^-)+\chi_r(\xi^+\},
$$
$$\mathcal{F}^- = 
\{F(\xi^-,\xi^+):u_x|_{[\xi^-,\xi^+]}=-1,\chi_l(\xi^-)+\chi_r(\xi^+\}.
$$
Since the number of elements in $\mathcal{F}^+ \cup \mathcal{F}^-$ is finite, we can order them,
$$F_1(\xi^-_1,\xi^+_1),\ldots, F_N(\xi^-_N,\xi^+_N),
$$
where
$$
a\le \xi^-_1 < \xi^+_1 \le \xi^-_2 < \xi^+_2 \le \ldots  \xi^-_i < \xi^+_i \le \ldots  \xi^-_N < \xi^+_N \le b.
$$
Facets occurring in this sequence can be grouped in the following way
$$F_i,\ldots, F_{i+l},$$
where $F_j\in \mathcal{F}^+,\ j=i,\ldots,i+l$ (or respectively, they belong to
$\mathcal{F}^-$) and $F_{i-1},\ F_{i+l+1}\in \mathcal{F}^-$ (respectively, they
are in $\mathcal{F}^+$) as far as $i>1,\ i+l<N$. 

Let us consider a typical group $F_i,\ldots, F_{i+l}$, for the sake of
definiteness, we assume that $F_j\in \mathcal{F}^+$, $j=i,\ldots, i+l$.

We prove our theorem by induction with respect to $l$. We notice  that $l$ is always odd.
Let $l=1$, We will estimate the  extinction time using the following procedure. We take 
$$
x_1 \in [\xi_i^-, \xi_i^+],\qquad x_2 \in [\xi_{i+1}^-, \xi_{i+1}^+],
$$
such that $u_x(x_1)=u_x(x_2)=1$ and we set
$$x_0 = \max\{y:u_x(y) = -1\ \wedge\ y<x_1\},\quad x_3 = \min\{y:u_x(y) = -1\ \wedge \ x_2<y \}.$$

It follows from the definition of $x_0$,  $x_3$ that the graph of $u_0$ restricted to interval $[x_0,x_3]$ is contained in the strip limited by the tangents at the points $x_0,\ldots,x_3$ (see Figure 2). We denote these tangent lines by $\ell_0,\ldots, \ell_3$, i.e. $(x_i, u_0(x_i))\in \ell_i$, $i=0,\ldots, 3$.
\begin{figure}[h]
 \begin{center}
{\epsfxsize=10cm \epsfbox{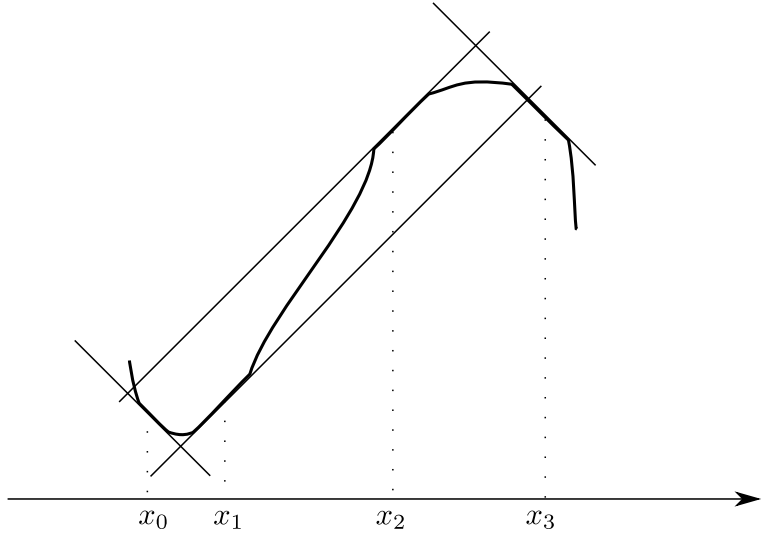}}

\end{center}
\caption{Case l=1.}
\end{figure}
First, we calculate the upper and lower estimates of function $u$.
We define $w\ge u \ge v$ in the following way:
$$w(x,\delta) = \left\{\begin{array}{cc}
				\min\{\ell_2,\ell_3\}, & x\in [x_2,x_3],\\
				u(x,\delta),& x\notin [x_2,x_3],
				\end{array}\right. $$
$$v(x,\delta) =   \left\{\begin{array}{cc}
				\max\{\ell_0,\ell_1\}, & x\in [x_0,x_1],\\
				u(x,\delta),& x\notin [x_0,x_1].
				\end{array}\right. $$
Consider solutions $w,\ v$ of (\ref{rng}), $t\geq \delta$ with initial conditions $w_0 = w(x,\delta), \ v_0 = v(x,\delta)$, respectively.
By the comparison principle, we have
$$v(x,t) \le u(x,t) \le w(x,t)\quad \textrm{for}\ t\ge \delta, \ x\in I.$$

We denote by $p_1(t)$, (resp. $p_2(t)$) the line parallel to $\ell_1$, , (resp. $\ell_2$) and passing through the point $(x_1,v(x_1,t)),$  (resp. $(x_1,w(x_1,t))$.
We can estimate the time $t_1$ such that $p_1(t_1) = p_2(t_2)$.
Indeed, if $\frac{dh_1}{dt},\ \frac{dh_2}{dt}$ are vertical velocities of $p_1,\ p_2$, respectively, then 

\begin{equation}\label{vv}
\int_\delta^{t_1}(\frac{dh_1}{dt}-\frac{dh_2}{dt} )= d,
\end{equation}
where $d=\ell_2(x_1)-\ell_1(x_1)$. We recall formula (\ref{rn-4.1}), 
$$\frac{dh_1}{dt} = \frac{2}{L_1},\quad \frac{dh_2}{dt} = -\frac{2}{L_2},$$
where $L_1$ is the length of facet at $(x_1,v(x_1,t))$, $L_2$ is the length of facet at $(x_2,w(x_2,t))$.

Note that $L_1+L_2 \le \lambda$, where $\lambda$ is the distance between projections of points $\ell_0\cap \ell_1$ and $\ell_2\cap \ell_3$ on the $x$-axis.
It follows form (\ref{vv}) that

$$d=\int_{\delta}^{t_1}(\frac{dh_1}{dt}-\frac{dh_2}{dt}) = \int_{\delta}^{t_1}(\frac{2}{L_1}+\frac{2}{L_2})  
\ge \int_{\delta}^{t_1}\frac{4}{\lambda}.$$

As a result,
$$
T(x_1,x_2) \le (t_1-\delta) \leq \frac{d\lambda}{4} = \frac{\lambda}{4}(\ell_2(x_1) - \ell_1(x_1).
$$
Next, we proceed inductively. Suppose that $l>1$. 
We denote by $\ell_j$ lines containing facets $F_j$, $j=i,\ldots,i+l$. We choose $k$ and $m$ such that

$$|\ell_k(x)-\ell_m(x)| = \max\left\{|\ell_r(x)-\ell_s(x)|:\ r,s \in\{i,\ldots,i+l\}\right\}.$$

This means that line $\ell_j$ containing $F_j$ lies within the strip bounded by $\ell_k,\ \ell_m$. 
There are two cases to consider:
\begin{itemize}
\item[(i)] there exists a facet $F_j$ contained within the strip bounded by $\ell_k,\ \ell_m$;
\item[(ii)] all other facets are contained in the $\ell_k$ or $\ell_m$. In such a case let $m$ be the largest and $k$ the smallest index that fulfills this assumption. 
\end{itemize}
\begin{figure}[h]
 \begin{center}
{\epsfxsize=7cm \epsfbox{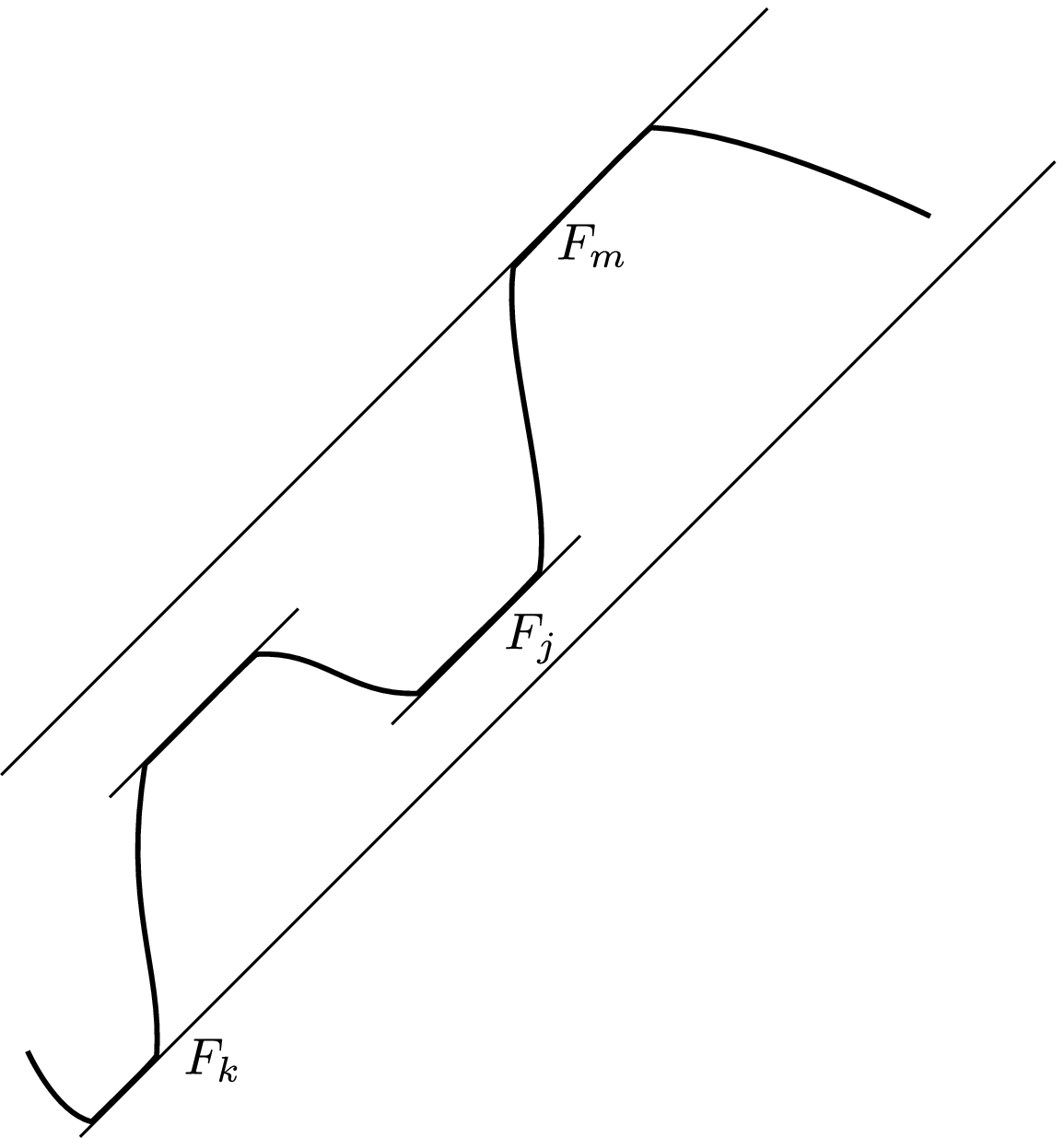}}

\end{center}
\caption{Case (i).}
\end{figure}
We first consider  case (i). 
Since $F_j = F(\xi^-_j,\xi^+_j)$ and for some $\epsilon >0$, $u|_{[\xi^-_j - \epsilon,\xi^+_j+\epsilon]}$ is on one of the side of $\ell_j$, it follows that in the group $F_i,\ldots, F_{i+l}$, except for $F_j,F_k,F_m$, there is one additional facet $F_{j+\epsilon}$ adjacent to $F_j$, where $\epsilon = 1$ or $\epsilon =-1$.
We are looking for a point of intersection of $\ell_j$ and the graph of $u$  
such that
$$\bar{x} = \max\{\tilde{x}<\xi^-_j:  \ell_j(\tilde{x})=u(\tilde{x},\delta)\}\quad \textrm{if}\ \epsilon =-1$$
or
$$\bar{x} = \min\{\tilde{x}>\xi^+_j:  \ell_j(\tilde{x})=u(\tilde{x},\delta)\}\quad \textrm{if}\ \epsilon =1.$$
We consider
$$v(x,\delta) = \left\{\begin{array}{cc}
				u(x,\delta),& x\notin [\bar{x},\xi^-_j],\\
				\ell_j(x),& x\in [\bar{x},\xi^-_j].
				\end{array}\right.$$
Then $v$ has $l-2$ facets with non-zero curvature. Function $w(x,\delta)$ is defined analogously, i.e.
we set, if $\epsilon =-1$
$$
\tilde x = \min\{ y> \xi^+_{j-1}:\ \ell_{j-1} (x) = u(x,\delta)
$$
or, if  $\epsilon =1$,
$$
\tilde x = \max\{ y> \xi^+_{j+1}:\ \ell_{j+1} (x) = u(x,\delta).
$$
For  $\epsilon=-1$, we define,
$$
w(x,\delta) = \left\{\begin{array}{cc}
				u(x,\delta),& x\notin [\xi^+_{j-1},\tilde x]\\
				\ell_j(x),& x\in [\bar{x},\xi^-_j].
				\end{array}\right.
$$
A similar definition is for 	 $\epsilon =1$.			
We have $$w(x,\delta)\ge u(x,\delta)\ge v(x,\delta).$$

In case (ii) we proceed analogously, except that  we now only have lines $\ell_k$ and $\ell_m$. (To define $v$ we combine the facet, which creates $\ell_k$, with the nearest facet with the same curvature). We see that $w$ and $v$ have $l-2$ facets with non-zero curvature. 
We shall use the induction hypothesis that we have the estimate of extinction
time for function $u$ in case $l-2$. Thus, we arrive  at an estimate for
$T(x_k,x_m)$, 
$$
T(x_k,x_m) \le \frac{\lambda}{4}(\ell_m(x_k) - \ell_k(x_k)),
$$
where $\lambda$ is the distance from the $x_1$-coordinate of $\ell_0\cap \ell_m$ to $\ell_m\cap \ell_{m+1}$, where $l_0$ is the line passing through $F_{i-1}$ and $\ell_{m+1}$ is the line including  $F_{i+l+1}$.

First, we will estimate the  time $T(x_k,x_m)$ after which $v$ and $w$ collide
in the strip bounded by lines $\ell_k$ and $\ell_m$. 
For this purpose we take
$$\underline{v}\le v \le \underline{w},\quad \overline{v}\le w \le \overline{w}.$$

We notice that our estimates on $T(x_k,x_m)$ 
we are
developing here depend only 
on the parameters of the
strip determined by the lines $\ell_k$,  $\ell_m$ bounding a part of the graph
of solution $u$. By construction, the same bounding box is for $u$ and $v$, $w$
yielding the same estimate for collision 
time 
$T(x_k,x_m)$. The estimate is made on the premise that lines 
$\ell_k$, $\ell_m$ sweep the strip. 

We note that in the case of Dirichlet 
boundary conditions if facet $F$ touches the boundary, then $F$ has zero curvature. We proceed as earlier
with one difference. In equation (\ref{vv}) we have only one nonzero vertical velocity. But this in some cases gives us even better extinction time. In the case of Neumann data we extend $u_0$ by odd reflection. This reflection does not lengthen the extinction time. Hence the claim. \qed

\begin{remark}
The inspection of the proof shows that it does not require differentiability of
solution. In our calculations we depended on the fact that for typical $t>0$,
all facets are created  and there are no missing directions
(see Theorem \ref{proB} and Proposition \ref{cc}). 
Thus, the argument remains valid also when the initial condition $u_0$ is in general in $BV$
for typical time $\delta>0$. 
\end{remark}

\subsection{Examples} \label{ssdisc}

Here, we present two examples highlighting the main issues addressed in this
paper.

\bigskip\noindent
{\bf Example 1.} In this example we study solutions when the initial data are discontinuous at $x=a$ or $x=b$. We see that the solution does not satisfy the boundary condition in a pointwise manner for $t>0$. 
%
\begin{proposition} \label{pr-ex}
Let us suppose that $I=(-1,1)$ and $u(t,-1) = u(t,1) = 0$.\\
(a) If $u_0(x) = -|x| +d$,
where $d>1$, then
$$
u(x,t) = -|x| +d -2\min\left\{t,\frac{d-1}{2}\right\}
$$
and $\Omega(u;x,t) = -2x.$
In particular, $T_{ext} = \frac{d-1}{2}$ and $u(x, T_{ext}) = 1 -|x|$.\\
(b) If $u_0(x) = -|x| + e$, where $e<1$, then
\begin{equation}\label{rn-kon}
 u(x,t) = \max\left\{u_0(x),|x|-2+e+2\min\left\{\sqrt{2t},\frac{(1-e)}{2}\right\} \right\}
\end{equation}
and $T_{ext} = \frac{(1-e)^2}{8}$.
In particular,  
$$u(x, T_{ext}) =\left\{\begin{array}{cc}
-x-1,& x\in [-1,-\frac{(1+e)}{2}],\\
                 -|x| + e, & x\in [-\frac{(1+e)}{2},\frac{(1+e)}{2}],\\
x-1,& x\in [\frac{(1+e)}{2},1].
                \end{array}\right.$$

\end{proposition}

\noindent{\it{Proof.}} We conduct calculations similar to that in the proof of Proposition \ref{cc}.
(a) We are interested in how long it takes for the facets to reach the extinction time. We have $h(0) = d$,  $h(T_{ext}) = 1$, the facets have constant length, $L(h) = 1$. We notice that $\Omega(x,t) = -2|x|$, hence
$$\frac{dh}{dt}(t) = \frac{-2}{L(h)}\quad i.e. \quad \frac{dh}{dt}(t)L(h) = -2.$$
We integrate the above equation over $[0,T_{ext}]$
$$\int_0^{T_{ext}}\frac{dh}{dt}(t)L(h)dt =-2T_{ext}.$$
We see that $T_{ext}= \frac{d-1}{2}$.

(b) We see that  due to Theorem \ref{proB} at points $x=-1$, and $x=1$ two symmetric facets are created, $F(-1,-\eta)$, $F(\eta, 1)$. We will follow $F(\eta, 1)$, hence 
we have $h(0)=-1 + e$, $h(T_{ext}) = 0$. The length of the facet is $L(h) =
\frac{1}{2}(1+h-e)$. At the same time $L= 1-\eta(t)$, as a result,
$$
\Omega (x,t) = 
\left\{\begin{array}{cc}
                 \frac{2x-2\eta(t)}{1-\eta(t)}, & x\in[\eta(t),1],\\
0,& x\in[-\eta(t),\eta(t)],\\
\frac{2x + 2\eta(t)}{1-\eta(t)},& x\in[-1,-\eta(t)]. 
                \end{array}\right. .
$$
Using this information we write equation for $h$,
$$\frac{dh}{dt}(t) = \frac{2}{L(h)}\quad i.e.\quad \frac{dh}{dt}(t)L(h) = 2.$$
We integrate the above equation to get 
$$
h(t) = 2 \sqrt{2t} - 1+e\quad\hbox{and}\quad \eta(t)= 1 - \sqrt{2t}. 
$$
Thus, we see that $T_{ext}=\frac{(1-e)^2}{8}$ and (\ref{rn-kon}) holds.\qed
\bigskip

It is worthwhile to consider an example of an oscillating initial condition.

\noindent
{\bf Example 2.} We look at a solution with oscillating initial data. Theorem \ref{T-ext} implies that facets close to zero get killed first, so that the extinction time is estimated by using the parameters corresponding to the biggest humps in the data.\\

We notice that $u_1(x) = x\sin(x^{-1})$ is not in $BV(I)$. On
the other hand, $u_0(x) = x^2\sin(x^{-1})\in BV(I)$. Let us suppose that $I=(-1,1)$. We see that for any $t>0$
most of the facet interaction are over, only a finite number of facets with
non-zero curvature are left. We approximate $u_0$ with
$$u_0^n(x) = \left\{\begin{array}{cc}
                  0,& x\in[-\frac{1}{n\pi},\frac{1}{n\pi}],\enskip\\
		  x^2\sin(x^{-1}),& x\in [-1,1]\setminus [-\frac{1}{n\pi},\frac{1}{n\pi}].
                 \end{array}\right.$$
Due to Proposition \ref{prop-Text} and Theorem \ref{T-ext} we have an estimate on
the extinction time for the evolution with initial condition $u_0$. However, we provide no closed formula for it.

\section{Viscosity solutions}\label{lepkie}
There are two main reasons for introducing the theory of viscosity solutions in this paper. Firstly, we would like to check if the  oscillatory behavior of solutions is \textquoteleft correct\textquoteright. At the same time the theory of viscosity solutions will give us an additional tool like the comparison principle, see Theorem \ref{comprin}, which is used in the proof of Theorem \ref{T-ext}. This is the second reason for dealing with the theory of viscosity solutions.

Our exposition is based on \cite{miyory}, adapted to the setting of (\ref{rng}), while avoiding unnecessary  generality of that paper. It is clear that we have to give meaning to $\left( {\cL}(\varphi_x)\right) _x$ for a proper choice of test functions $\varphi$. Our experience with the theory of nonlinear semigroups suggests that it is advantageous to work with $\left( W_p(\varphi_x)\right) _x$, when $W(p)$ is a convex function. 

We are going to present the necessary notions. 
A function $f\in C(I)$ is called {\it faceted} at $x_{0}$ with 
{\it slope} $p\in\{-1, 1\}=:\cP$ on $I$ (or $p$-{\it faceted} at $x_{0}$) 
if there is a closed nontrivial finite interval $\tilde I \subset I)$ containing 
$x_{0}$ such that $f$ coincides with an affine function 
$$
\ell_{p}(x)=p(x-x_{0})+f(x_{0}) \ {\rm in} \ \tilde I
$$
and $f(x)\neq \ell_{p}(x)$ for all $x\in J \backslash \tilde{I}$, where $J$ is a
neighborhood of $\tilde I$ in $I$. 
The interval $\tilde I$ 
is denoted by $R(f, x_{0})$. 

We will denote by  $C^{2}_{P}(I)$ the set of $f\in C^{2}(I)$ such that $f$ is 
$p$-faceted at $x_{0}$ whenever $f^{\prime}(x_{0})\in \{-1, 1\}$.
Let $A_{P}(I_T)$ be the set of all admissible functions $\psi$ on $I_T$ i.e. $\psi$ is of the form
$$\psi(x,t) = f(x)+g(t),\quad f\in C^{2}_{P}(I),\ g\in C^1(0,T).$$

The definition of $\left( W_p(\varphi_x)\right) _x$ is non-local for $p$-faceted $\varphi\in C^{2}_P(I)$. It involves a solution of an obstacle problem, which we will describe momentarily. 

We assume that 
$\Delta>0$, $\chi_l$, $\chi_r\in\{1,-1\}$ and 
$J= [\alpha,\beta]\subset I$ are given. We set
\begin{eqnarray*}
K^Z_{\chi_l\chi_r}(J) &=&
\left\{ \zeta\in H^1(J): \ Z(x)-\frac{\Delta}{2} \le \zeta(x) \le Z(x) +\frac{\Delta}{2}, \ x\in J, \right. \\ && \left.
Z(\alpha)- \chi_l\frac{\Delta}{2} = \zeta(\alpha), \ Z(\beta)+ \chi_r \frac{\Delta}{2}= \zeta(\beta)\right\}.
\end{eqnarray*}
We also introduce
$$
\cJ^Z_{\chi_l\chi_r}(\zeta,J) =
\begin{cases}
 \int_J |\zeta'(x)|^2 \,dx & \hbox{if }\zeta \in H^1(J),\\
+\infty & \hbox{if }\zeta \in L^2(J)\setminus H^1(J).
\end{cases}
$$
Let us call by $\xi^{Z,J}_{\chi_l\chi_r}$ the unique solution to the obstacle problem
\begin{equation}\label{min-p}
 \min\{ \cJ^Z_{\chi_l\chi_r}(\zeta,J) : \ \zeta\in K^Z_{\chi_l\chi_r}(J)\}.
\end{equation}
It is easy to see that for any affine function $Z$, the minimizer is
an affine function too. 
This is the case considered in this paper. But in  general, even if $Z\in C^2$, then it is
well-known that the unique solution $\xi^{Z,J}_{\chi_l\chi_r}$ belongs
to $C^{1,1}(J)$, see \cite{kinder}.

Of course, $Z$ is defined up to an additive constant, which we
have to choose  properly. 
If $F=F(\xi^-,\xi^+)$ is a facet, and 
$u_x\left|_{(\xi^-,\xi^+)} \right. =p\in \cP$, then we set $Z=p x$ and  
$\Delta = W_p(p^+) - W_p(p^-)=2$. 

We define $\Lambda^Z_W(\varphi )$  as follows. We stress that in
\cite{miyory} we denoted 
the same object by $\Lambda^{Z'}_W(\varphi)$, but here we opt for a simpler notation.
 
If $\varphi\in C^2$ and $\varphi_x\not\in\cP$, then we set
$$
\Lambda^Z_W(\varphi )(x) := \left( W_p(\varphi_x)\right) _x .
$$
If $\varphi\in C^2_P$ is $p$-faceted at $x_0$, then we denote its
faceted region $R(\varphi,x_0)$ by $J$. We take $Z := p\in \cP$, then we set
$$
\Lambda^Z_W(\varphi )(x) := \frac{d}{dx} \xi^{Z,J}_{\chi_l\chi_r}(x).
$$
Transition numbers $\chi_l,\ \chi_r$ are defined by (\ref{tr-n1}) and (\ref{tr-n2}).

It turns out that the non-local definition of $\Lambda^{Z}_W(\varphi )$ has the desired property.

\begin{proposition}(\cite[Theorem 2.4]{miyory})\label{comp}
 Assume that $I_{1}$ and $I_{2}$ are bounded open intervals and $\xi^{Z,I_i}_{\chi_l\chi_r}$, $i=1,2$ is the solution to  (\ref{min-p}).  We write $\Lambda_{\chi_l\chi_r}(x,J)$ for 
$\frac{d}{dx} \xi^{Z,J}_{\chi_l\chi_r}(x)$.
\begin{itemize}
\item[(i)]\quad
If $I_2 \subset I_1$, then
\begin{equation*}
\Lambda_{--}(x,I_2) \leq \Lambda_{\pm\pm}(x,I_{1}) \leq \Lambda_{++}(x,I_{2}) \
{\rm for} \ x\in I_{2}.  
\end{equation*}
\item[(ii)]\quad
If $a \leq c < b \leq d$ for $I_1 =(a,b)$, $I_2=(c,d)$, then for $x \in (c,b)$ 
\begin{equation*}
 \Lambda_{\pm-}(x,I_1) \leq \Lambda_{+\pm}(x,I_2), \enskip
\Lambda_{-\pm}(x,I_2)\leq \Lambda_{\pm+}(x,I_1). 
\end{equation*}
\end{itemize} 
\end{proposition}

After these preparations we may define the test functions and viscosity solutions to (\ref{rng}).

\begin{definition}{\rm
A real-valued function $u$ on $I_T$ is a (viscosity) {\it subsolution} 
of (\ref{rng}) in $I_T$ if
its upper-semicontinuous envelope $u^{*}$ is finite in  $\bar I_T$ and 
\begin{equation}\label{3.3}
\psi_{t}(\hat{t},\hat{x}) -  
\Lambda^{Z(\hat{t}, \cdot)}_{W} (\psi(\hat{t})) \ (\hat{x})  \leq 0 
\end{equation}
whenever $\left( \psi,(\hat{t},\hat{x}) \right) \in A_{P}(I_T)\times I_T$ fulfills
\begin{equation}\label{3.4}
\max_{I_T} (u^{*}-\psi)=(u^{*}-\psi) \ (\hat{t},\hat{x}).  
\end{equation}
Here{,} $\psi(\hat{t})$ is a function on $\Omega$ defined by
$\psi(\hat{t})=\psi(\hat{t}, \cdot)$
and $u^{*}$ is defined by
$$
u^{*}(t,x)=\lim_{\varepsilon \downarrow 0} \ \sup \{ u(s,y) | \ |s-t|<\varepsilon, \
|x-y|<\varepsilon, \ (s,y)\in I_T \}\quad\hbox{for }(t,x)\in\overline{I_T}.
$$
We also set $u_{*}=(-u^{*})$.

A (viscosity) {\it supersolution} is defined by replacing $u^{*}(<\infty)$ by the lower-semicontinuous envelope
$u_{*}(>-\infty)$, max by min in (\ref{3.4}) and the inequality (\ref{3.3}) 
by the opposite one. 
If $u$ is both a sub- and supersolution, it
is called a {\it viscosity solution}. 
Hereafter, we avoid using the word viscosity, if there is no ambiguity. 
Function $\psi$ satisfying (\ref{3.4}) is called a test function of $u$ at
$(\hat{t}, \hat{x})$. }
\end{definition}

The main tool which we acquire from the theory of viscosity solutions is the Comparison Principle.
\begin{theorem}(\cite[Theorem 4.1]{miyory})\label{comprin}
Let $u$ and $v$ be respectively sub- and supersolutions of (\ref{rng}) 
in $I_T=I \times(0,T)$, where $I$ is a bounded open interval. 
If $u^{*}\leq v_{*}$ on the parabolic boundary 
$\partial_{p}I_T(=[0,T)\times\partial I \cup \{ 0\}\times \overline{I})$ of $I_T$, 
then $u^{*} \leq v_{*}$ in $I_T$.
\end{theorem}

In order to use Theorem \ref{comprin}, we have to check that the solutions we constructed in Theorem \ref{tw-gl}, are actually viscosity solutions.

\begin{theorem}\label{pr-nav}
 If $u_0\in BV$ and $u$ is the corresponding solution to (\ref{rng}) with either boundary condition (\ref{rnD}), (\ref{rnP}) or (\ref{rnN}). Then $u$ is a viscosity solution to (\ref{rng}).
\end{theorem}

Before we  engage into the proof of Theorem \ref{pr-nav}, we make observations
facilitating the argument. We set,
$$
\Xi(u) =\{ [\xi^-,\xi^+]\subset [a,b]: [\xi^-,\xi^+]\hbox{ is the
pre-image of facet }F(\xi^-,\xi^+)\hbox{ of }u \}.
$$

\begin{lemma}\label{lm-a}{\rm 
 Let us suppose that $u(\cdot,t_0)$ is continuous at $x_0$ and $u_t$
 exists at $(x_0,t_0)$, where $t_0>0$. Then one of the following
 possibilities holds:\\  
(a) $x_0$ is in the complement of the sum of all pre-images of facets, i.e.
$x_0 \in I \setminus \bigcup \Xi(u(\cdot,t))$. Hence, $u_t(x_0,t_0) = 0$.\\ 
(b) $x_0$ is in the interior of the  pre-image of a facet, $x_0\in
 (\xi^-(t_0),\xi^+(t_0))$ and either (i)  the length of the interval 
$\xi^+(t) -\xi^-(t)$  as a function of time is continuous at $t_0$ and 
$\chi_r +\chi_l$ is a non-zero constant for all $t$ from a
neighborhood of $t_0$ or (ii) $\chi_r +\chi_l =0$  for all $t$ from a
neighborhood of $t_0$. \\
(c) $x_0\in \{\xi^-(t_0),\xi^+(t_0)\} $ and either (i) facets
 $F(\xi^-(t),\xi^+(t))$ have zero curvature for all $t<t_0$
 sufficiently close to $t_0$
 or  (ii) $\chi_r +\chi_l \not =0$ and $x_0 \in
 F(\xi_1^-(t_0),\xi_1^+(t_0)) \cap F(\xi_2^-(t_0),\xi_2^+(t_0))$ and
 functions $\xi_1^+(\cdot) - \xi_1^-(\cdot)$, $\xi_2^+(\cdot) -
 \xi_2^-(\cdot)$ are continuous at $t=t_0$ and equal at $t=t_0$. We
 note that (i)  includes the case of a facet passing
 through the boundary data.}
\end{lemma}

\noindent{\it Proof.} If $u_t(x_0,t_0)$ exists, then this means that
$u^+_t(x_0,t_0) = u^-_t(x_0,t_0)$. We know how to calculate
$u^+_t(x_0,t_0)$. First, we consider the case of $x_0$ belonging to
$I\setminus \bigcup \Xi(u(\cdot,t_0))$. Since facets with non-zero curvature
are
expanding, then   $x_0 \in I\setminus \bigcup \Xi(u(\cdot,t))$
for all $t<t_0$. Since   $u_t(x_0,t_0)$ exists, we deduce that it must
be zero. Thus, (a) holds.

Let us now assume that $x_0$ belongs to the interior of the 
pre-image of a facet, $x_0\in(\xi^-(t_0), \xi^+(t_0))$. Then always, 
\begin{equation}\label{rn-p}
u^+_t(x_0,t_0) = \frac{\chi_r +\chi_l}{\xi^+(t_0)- \xi^-(t_0)}.
\end{equation}
Continuity of $t\mapsto  u(x_0,t)$ at $ t_0$ implies that (b) holds.

Let us suppose that not only $x_0\in\{\xi^-,\xi^+\}$ but also
$x_0 \in \partial \left(\bigcup \Xi(u(\cdot,t_0))\right)$. If $\chi_l(\xi^-(t))+ \chi_r(\xi^+(t))\neq 0$
for for $t<t_0$ sufficiently close to $t_0$, then $x_0 \in
I \setminus \bigcup \Xi(u(\cdot,t))$, in this case, 
unless $\chi_r +\chi_l =0$, i.e. facet
$F(\xi^-, \xi^+)$ has zero curvature, including the case of a facet
satisfying the boundary conditions. Hence, $u_t(x_0,t_0) =0$.

If $(x_0, u(x_0,t_0)) \in
 F(\xi_1^-(t_0),\xi_1^+(t_0)) \cap F(\xi_2^-(t_0),\xi_2^+(t_0))$, then
 $u^+_t(x_0,t_0)$ does not exist unless 
$\xi_1^+(\cdot) - \xi_1^-(\cdot) = \xi_2^+(\cdot) -
 \xi_2^-(\cdot)$ and the transition numbers of 
$F(\xi_1^-(t_0),\xi_1^+(t_0))$ and $F(\xi_2^-(t_0),\xi_2^+(t_0))$ are
 the same. The existence of $u_t(x_0,t_0)$ and (\ref{rn-p}) imply the continuity of functions $\xi_1^+(\cdot) - \xi_1^-(\cdot)$, $\xi_2^+(\cdot) -
 \xi_2^-(\cdot)$  at $t=t_0$. Thus, we showed that (c) is satisfied.
\qed

We notice that fulfilling conditions (a) to (c) is also sufficient for the existence
of $u_t(x_0,t_0)$.

Equally important as Lemma \ref{lm-a} is understanding the
behavior of $u$ near $(x_0,t_0)$, when $u(x_0,\cdot)$ is not
differentiable with respect to $t$ at $t_0$.

\begin{lemma}\label{lm-b}{\rm
 Let us suppose that $u(\cdot,t_0)$ is continuous at $x_0$ and $u_t$
 does not exist at $(x_0,t_0)$. Then one of the following
 possibilities holds:\\ 
(a) $x_0$ is in the interior of a  pre-image of a facet, $x_0\in
 (\xi^-(t_0),\xi^+(t_0))$ and $\chi_l(t_0) +\chi_r(t_0)=0$, while
 $\chi_l(t) +\chi_r(t)\not=0$ for $t<t_0$ close to $t_0$. This case
 includes the situation when a facet hits the boundary of $I$ at time $t_0$.\\ 
(b) $x_0$ is in the interior of a  pre-image of a facet, $x_0\in
 (\xi^-(t_0),\xi^+(t_0))$ and $\chi_l(t) +\chi_r(t)=const\neq 0$, for $t\le t_0$
close to $t_0$ and $\xi^+(t)-\xi^-(t)$ has a jump at $t=t_0$.\\
(c) $x_0$ belongs to the boundary of $\bigcup \Xi(u(\cdot,t))$,
while  for all $t<t_0$ sufficiently close to $t_0$, point $x_0$
belongs to $I \setminus \bigcup \Xi(u(\cdot,t))$. In this case
$u^+_t(x_0,t_0)$ is given by  (\ref{rn-p}) and $u^-_t(x_0,t_0)=0$.\\
(d) $x_0 \in
 F(\xi_1^-(t_0),\xi_1^+(t_0)) \cap F(\xi_2^-(t_0),\xi_2^+(t_0))$ and
 functions 
$\xi_1^+(t_0) - \xi_1^-(t_0) \neq \xi_2^+(t_0) - \xi_2^-(t_0)$.
}
\end{lemma}

\noindent{\it Proof.} Here we listed cases complementary to those enumerated
in Lemma \ref{lm-a}, thus not further argument is necessary. 
\qed

We notice that  in case (d) function $x\mapsto u_t^+(x, t_0)$ is
discontinuous at $x=x_0$.

\bigskip
We are now ready for the 
{\it proof of Theorem \ref{pr-nav}.} We will show that $u$ is a subsolution. The argument that $u$ is also a supersolution is similar and it will be omitted. Let us take $(x_0,t_0)\in I_T$ and  $u^*$ the upper semicontinuous envelope of $u$. There are the following four cases to consider:\\
1) $x_0$ is a point of continuity of $u(\cdot,t)$, i.e. $u^*(x_0,t_0) = u(x_0,t_0)$;\\
2) $x_0$ is a point of discontinuity of  $u(\cdot,t)$. We note that jumps are the only possible discontinuities of $BV$ functions.

In each of the above situations either:\\
a) $x_0$ belongs to a pre-image of facet $F$, i.e. $x_0\in [\xi^-,\xi^+]$ or\\
b) the converse holds. 

\bigskip
{\it Step 1.}
We begin with  case 1 a). 
Let us
take a test function $\psi(x,t) = f(x) + g(t)$ such that $\psi(x_0,t_0) = 
u(x_0,t_0)$ and 
\begin{equation}\label{warlep}
\max (u -\psi) = u(x_0,t_0) - \psi(x_0,t_0). 
\end{equation}
We have to proceed according to the properties of the test
function. Let us assume first that $u_t(x_0,t_0)$ exists. Then
$g'(t_0) = u_t(x_0,t_0)$ and we need to consider the following cases. 


If $F$ has zero curvature with slope $p=1$ (case $p=-1$ is analogous), then we have the following cases:
\begin{itemize}
 \item[(i)] $F(R(f,x_0))$ has a nonzero curvature and $\xi^-\leq x^- < \xi^+ \leq x^+$ (or $x^-\leq \xi^-<x^+\leq \xi^+$), where $[x^-,x^+]=R(f,x_0)$. Then Proposition \ref{comp} implies $\Lambda_{+-}(x,[\xi^-,\xi^+]) \le \Lambda_{++}(x,R(f,x_0))$ (or $\Lambda_{-+}(x,[\xi^-,\xi^+]) \le \Lambda_{++}(x,R(f,x_0))$). Hence, (\ref{3.3}) holds.
\item[(ii)] $F(R(f,x_0))$ has zero curvature and $\xi^-\leq x^- < \xi^+ \leq x^+$ (or $x^-\leq \xi^-<x^+\leq \xi^+$). Then Proposition \ref{comp} implies $\Lambda_{+-}(x,[\xi^-,\xi^+]) \le \Lambda_{+-}(x,R(f,x_0))$ (or $\Lambda_{-+}(x,[\xi^-,\xi^+]) \le \Lambda_{++}(x,R(f,x_0))$). Hence, (\ref{3.3}) holds.
\item[(iii)] $F(R(f,x_0))$ has a nonzero curvature and $R(f,x_0)\subset [\xi^-,\xi^+]$. Then Proposition \ref{comp} implies $\Lambda_{+-}(x,[\xi^-,\xi^+]) \le \Lambda_{++}(x,R(f,x_0))$. Hence, (\ref{3.3}) holds.
\end{itemize}

Let us consider $x_0\in [\xi^-,\xi^+]$ and 
facet $F(\xi^-,\xi^+)$ has $\chi_l+\chi_r =-2$. In this case
$\Omega_x |_{(\xi^-,\xi^+)}= -2/(\xi^+-\xi^-)$ and $\Omega_x =
\Lambda_{--}(x,(\xi^-,\xi^+))$. If the graphs of $u$ and $f$ are below the line
$l_p$ passing through $F(\xi^-,\xi^+)$, then the faceted region $R(f,x_0)$
contains $[\xi^-,\xi^+]$. Then we consider $\Lambda_{--}(x,R(f,x_0))$. We also
notice that
$$
\Lambda_{--}(x,R(f,x_0)) \ge \Lambda_{--}(x,(\xi^-,\xi^+)).
$$
Hence, (\ref{3.3}) holds.

The other possible tangency configurations of $u$ and $f$ are analyzed as in
(i)--(iii) above. The details are left to the interested reader.

The case  $x_0\in [\xi^-,\xi^+]$ when 
$\chi_l+\chi_r =2$ is simpler, because there is
just one way how $f$ may touch $u$. Namely, we have $\chi_l=\chi_r=1$ and
$R(f,x_0)\subset [\xi^-,\xi^+]$. Thus, $\Lambda_{++}(x,R(f,x_0)) \ge
\Lambda_{++}(x,(\xi^-,\xi^+))$ and  (\ref{3.3}) follows.
\bigskip

{\it Step 2.} Now, we work assuming that $u_t(x_0,t_0)$ does not exist. We have two major
subcases: 
\begin{equation}\label{case-i}
x_0\in (\xi^-, \xi^+),
\end{equation}
\begin{equation}\label{case-ii}
x_0\in \{\xi^-,\xi^+\}.
\end{equation}
Our analysis is based on Lemma \ref{lm-b}. 
If (\ref{case-i}) and the case Lemma \ref{lm-b} (a) hold, then
$F(\xi^-, \xi^+)$ has zero curvature and there exist facets
$F(\zeta^-(t),\zeta^+(t))$, such that
\begin{equation}\label{zetapm}
 \lim_{t\to t_0^-}\zeta^\pm(t) =: \zeta^\pm \in [\xi^-,\xi^+].
\end{equation}
First, we consider
$$
x_0\in (\xi^-,\xi^+)\setminus \{\zeta^-,\zeta^+\}.
$$ 
We will separately study
\begin{equation}\label{left}
x_0\in \{\zeta^-,\zeta^+\}. 
\end{equation}
We have two obvious possibilities for $F(\zeta^-(t),\zeta^+(t))$, either
$\chi_l(t)+\chi_r(t) < 0$ or $\chi_l(t)+\chi_r(t) > 0$, here we use the
shorthands, $\chi_l(t) \equiv \chi_l(\zeta^-(t))$, $\chi_r(t) \equiv
\chi_l(\zeta^+(t))$.

If $\chi_l(t)+\chi_r(t) < 0$ for $t<t_0$ close to $t_0$, then there is
no $g\in C^1(0,T)$ such that
\begin{equation}\label{rn-c}
f(x_0)+ g(t_0) = u(x_0,t_0),\quad u(x,t) \le g(t) + f(x)\quad
\hbox{in a neighbourhood of }(x_0,t_0).
\end{equation}
On the other hand, if 
$\chi_l(t)+\chi_r(t)>0$ for $t<t_0$ close to $t_0$, then there exists
$g$ satisfying (\ref{rn-c}) with
$g'(t_0)\in [0,A]$. 

If $x_0 \in [\zeta^-, \zeta^+]$, when $\zeta^\pm$ are defined by
(\ref{zetapm}), then by
Lemma \ref{lm-b}, we deduce that $A =
\frac{2}{\zeta^+ -\zeta^-}$.
We have further subcases to consider:\\
($\alpha$) $\zeta^+ =\xi^+$, $\zeta^- =\xi^-$, i.e. facet
$F(\xi^-,\xi^+)$ is a result of a collision of a facet moving upward with
the boundary of $I$.\\
($\beta$) facet
$F(\xi^-,\xi^+)$ is a result of a collision of a facet moving upward
with a facet passing through the boundary data.\\
($\gamma$) facet
$F(\xi^-,\xi^+)$ is a result of a collision of a facet moving upward
with another facet  moving downward.

Let us consider the resulting limitations on $f$ and $R(f,x_0)$. Case
($\alpha$) does not bring any. If ($\beta$) occurs, then (\ref{rn-c}) 
implies that either 
\begin{equation}\label{case-a}
 R(f,x_0)\subset [\zeta^-,\zeta^+]
\end{equation}
 or
\begin{equation}\label{case-b}
R(f,x_0)\not\subset [\zeta^-,\zeta^+] 
\end{equation}
but $g'(t_0)=0$. Finally,
($\gamma$) and  (\ref{rn-c})  imply that $R(f,x_0)\subset
[\zeta^-,\zeta^+]$, because the situation is similar to that studied in the
lines above formula (\ref{rn-c}).

Let us check that $u$ is a subsolution in these cases. We notice that
\begin{equation}\label{rn-d}
\Omega_x (x_0,t_0)= 0 = u_t^+ (x_0,t_0).
\end{equation}
It will be easier if we start with ($\gamma$) first. If this case we have
$R(f,x_0)\subset [\zeta^-,\zeta^+]$ and
$\Lambda^Z_W(f,x) = \Lambda_{++}(x, R(f,x_0))$. Since $g'(t_0) \le
\frac{2}{\zeta^+ -\zeta^-} = \Lambda_{++}(x,[\zeta^-,\zeta^+]) $, then
we infer from Proposition \ref{comp} that $\Lambda_{++}(x, R(f,x_0)) \ge
\Lambda_{++}(x,[\zeta^-,\zeta^+]) $, thus (\ref{3.3}) holds.

In case
($\alpha$) there is no apparent restriction on $R(f,x_0)$ but it is
not clear, which minimization problem is the correct one if
$R(f,x_0)$ intersects the boundary of $I$. Since we developed ideas
for the Dirichlet boundary condition through the periodic boundary
data, we first extend $u$ in the way we did it earlier to get a periodic
function, cf. (\ref{war-ok}). We see that ($\alpha$) corresponds to ($\gamma$)
considered above.
Thus, we immediately conclude that if ($\alpha$) holds, then (\ref{3.3}) is
satisfied as well. We also check it directly. 
Formula (\ref{war-ok}) implies
that
($\alpha$) corresponds to
the collision of two facets, one is moving
downward the other is moving upward. In this case, $\Lambda^Z_W(f,x) =
\Lambda_{++}(x, R(f,x_0))$. Since we have
$$
\Omega_x (x,t_0)= \Lambda^Z_W (u,x)=  0, 
$$
and $g'(t_0) \le
\frac{2}{\zeta^+ -\zeta^-} = \Lambda_{++}(x,[\zeta^-,\zeta^+])$, then
we infer from Proposition \ref{comp} that 
$$\Lambda_{++}(x, R(f,x_0)) \ge
\Lambda_{++}(x,[\zeta^-,\zeta^+]), $$ thus (\ref{3.3}) indeed holds.

We turn our attention to ($\beta$).
If the subcase (\ref{case-a}) holds, then $R(f,x_0)\subset [\zeta^-,\zeta^+]$ and
$\Lambda^Z_W(f,x) = \Lambda_{++}(x, R(f,x_0))$. Since $g'(t_0) \le
\frac{2}{\zeta^+ -\zeta^-} = \Lambda_{++}(x,[\zeta^-,\zeta^+]) $, then
we infer from Proposition \ref{comp} that 
$$
\Lambda_{++}(x, R(f,x_0)) \ge
\Lambda_{++}(x,[\zeta^-,\zeta^+]) ,$$
thus (\ref{3.3})
holds.

In subcase (\ref{case-b}), $R(f,x_0)\subset [\xi^-,\xi^+]$ and  if $R(f,x_0)$
does not intersect
$\partial I$, then $\Lambda^Z_W(f) = \Lambda_{++} (x,R(f,x_0))\ge \Lambda^Z_W(u).$
On the other hand, 
if $R(f,x_0)$ intersects $\partial I$, then we proceed as above in case
($\alpha$) and
take $\Lambda_{++} (x,R(f,x_0))$ for $\Lambda^Z_W(f)$. Hence,
$\Lambda^Z_W(f)\ge \Lambda^Z_W(u)$. As a result, in both cases (\ref{3.3})
holds. 

Now, we come back to the left out case (\ref{left}). If  
$F(\xi^-,\xi^+)$ does not
touch the boundary, then for $t<t_0$ close to $t_0$, $x_0$ does not belong to
any pre-image of any facet. This is so because $\zeta^\pm(t)$ are not constant,
see (\ref{warXI}), unless $\zeta^\pm(t)$ are points of discontinuity
of $u$. 
Hence,  
there is
no test function satisfying (\ref{warlep}). 

The other case is that $F(\zeta^-,\zeta^+)$ intersects the boundary.
As usually,
we have two possibilities for this facet, either
$\chi_l(t)+\chi_r(t)<0$ or $\chi_l(t)+\chi_r(t)>0$.
If $\chi_l(t)+\chi_r(t)<0$, then there is no $g\in C^1(0,T)$ such that
(\ref{warlep}) holds. If $\chi_l(t)+\chi_r(t)>0$, then we proceed as
in previous cases.
\bigskip

{\it Step 3.}
We consider the situation when  (\ref{case-i}) and the case Lemma \ref{lm-b} (b) hold.
Thus, a moving facet collides with a zero curvature facet. We have the
situation similar to that in Step 1. Thus, we may rule out the case of
$(\chi_l + \chi_r)(t^-) =- 2$ as impossible to satisfy (\ref{rn-c}).

If $(\chi_l + \chi_r)(t^+) = 2$, we conclude that the only possibility for
$R(f,x_0)$ is that  $R(f,x_0)\subset [\zeta^-,\zeta^+]$ and we have
$\Lambda^Z_W(f) = \Lambda_{++}(x,R(f,x_0))$. Arguing as before we see that
$g'(t)\in [0,2/(\xi^+-\xi^-)]$ and
$2/(\xi^+-\xi^-)]=\Lambda_{++}(x,[\xi^-,\xi^+]),$ but $\Omega_x=0$. As a result,
(\ref{3.3})
holds.

\bigskip

{\it Step 4.} Let us assume that (\ref{case-ii}) and  case (c) of Lemma \ref{lm-b}
hold. But there is no test function $\psi(x,t) = f(x) + g(t)$ such that
(\ref{rn-c}) holds. 

\bigskip

{\it Step 5.} Let us assume that (\ref{case-ii}) and  case (d) of Lemma \ref{lm-b}
hold.
If $F(\xi^-,\xi^+)$ has positive curvature, i.e. $\chi_l+\chi_r>0$, then Lemma
\ref{lm-b} (d) and
(\ref{rn-c}) imply that there is no test function.
On the other hand, i.e. if $\chi_l+\chi_r<0$, %
then there are test functions. In
this configuration $u^+_t(x_0,t_0) <0$ and $u^-_t(x_0,t_0) =0$. 
We may assume that $(\xi,u(\xi))$ is a common point of two facets
$F[\xi^-,\xi]$ and $F[\xi,\xi^+]$ and without the loss of
generality, 
$F(\xi^-,\xi)$ has zero curvature while
$F(\xi,\xi^+)$ has negative curvature, i.e. $\chi_l+\chi_r =-2$.

We deduce, that
$g'(t_0)\in [A,0]$, where
$$
A= \frac{-2}{\xi^+ - \xi}.
$$
Moreover, $f$ may be faceted with facets $p=\pm 1$ as well as $|f'(x_0)|<1$.
Then, it is easy to check that  (\ref{3.3}) holds.

Due to Lemma \ref{lm-a} and \ref{lm-b} all cases corresponding to 1a) are
exhausted.




\bigskip

{\it Step 6.}
Let us now consider situation corresponding to 2a) and its consequences.
If this occurs, then $x_0$ belongs to a facet $F=F(\xi^-,\xi^+)$ and $u$ is
discontinuous at $x_0$. This discontinuity implies that $x_0$ as an endpoint of
facet $F(\xi^-(t),\xi^+(t))$ does not move for $t$ in a neighborhood of $t_0$.

Furthermore, it may happen that
$u_t(x_0,t_0)$ exists. Then our argument is similar to that used
in Step 1, while taking into account that $|\Omega(x_0,t_0)|=2$ and
$u^*(x_0,t_0) = u(x_0,t_0)$ or $u^*(x_0,t_0) \neq u(x_0,t_0)$. The details are
left to the interested reader except for a new situation arising when the test
function has slope different from $\pm 1$. For the sake of definiteness we
assume that the slope of $F$ is 1. If $u_t(\cdot,t)$ is additionally continuous
at $x_0$, then $ u_t(x_0,t_0)=0$. If $u_t(\cdot,t)$ is discontinuous at $x_0$,
then $ u_t(x_0,t_0)= - 2/(\xi^+ - \xi^-)$. Any non-faceted test function
$\psi(x,t) = f(x) + g(t)$ must be such that $f'(x_0)> 1$. If this happen that
$\Lambda_W^Z = (W_p(f'(x)))_x|_{x=x_0} = 0$. Hence, (\ref{3.3}) holds.

The case when $F$ has slope $-1$ is handled in the same way.

If $u_t(x_0,t_0)$ does not exist, then we have several sub-cases:
\begin{itemize}
 \item[${1^*}$] facet is an effect of the collision of $F(x_0,\xi^+)$ with $F(\zeta^-,\zeta^+)$. Furthermore, $F(\zeta^-,\zeta^+)$ may have positive or zero-curvature;
\item[${2^*}$] $F(a,x_0)$ is an effect of collision of $F(\zeta^-,x_0)$, $(p=-1)$, with the boundary;
\item[${3^*}$] $F(a,x_0)$ is an effect of collision of $F(\zeta^-,x_0)$, $(p=-1)$, with a facet touching the boundary.
\end{itemize}
Those situations are analogous to that considered in Step 2, where we have
discontinuity of $u$ at $t_0$. We also  note that discontinuity of $u$ may lead
to non-faceted test functions as in the previous paragraph. The details,
however, are left to the interested reader.


The cases 1b) and 2b) are now easy and they are left to the reader. 

\qed

\section*{Acknowledgement}
A part of the research was performed during the visits of MM to the University of Warsaw supported by the
Warsaw Center of Mathematics and Computer Science through the program `Guests'.
Both authors enjoyed partial support of the NCN through 2011/01/B/ST1/01197 grant.

\end{document}